\documentclass[12pt]{article}
\usepackage[leqno]{amsmath}
\usepackage{amssymb,mathrsfs,bm}
\usepackage{hyperref}
\usepackage{mathtools}
\usepackage{fullpage}
\usepackage{epic}
\usepackage{eepic}
\usepackage[english]{babel}
\usepackage[all]{xy}
\usepackage{enumerate}
\usepackage[T1]{fontenc}
\usepackage[latin1]{inputenc}
\usepackage{color}
\usepackage{soul}
\usepackage{mathdots}
\usepackage{mathabx}
\usepackage{dsfont}
\newtheorem{theointro}{Theorem}

\newtheorem{theo}{Theorem}[subsection]
\newtheorem{lem}[theo]{Lemma}

\newtheorem{prop}[theo]{Proposition}
\newtheorem{cor}[theo]{Corollary}
\newtheorem{rem}[theo]{Remark}

\def\R{\mathbb{R}}
\def\C{\mathbb{C}}

\DeclareMathOperator{\Hom}{Hom}

\DeclareMathOperator{\Tra}{Tr}

\DeclareMathOperator{\Temp}{Temp}

\DeclareMathOperator{\Ker}{Ker}

\DeclareMathOperator{\Ad}{Ad}

\DeclareMathOperator{\vol}{vol}
\DeclareMathOperator{\Supp}{Supp}

\DeclareMathOperator{\cusp}{cusp}

\DeclareMathOperator{\Irr}{Irr}

\DeclareMathOperator{\cS}{\mathcal{S}}

\DeclareMathOperator{\fS}{\mathfrak{S}}
\DeclareMathOperator{\cA}{\mathcal{A}}
\DeclareMathOperator{\bR}{\mathbb{R}}
\DeclareMathOperator{\bZ}{\mathbb{Z}}
\DeclareMathOperator{\bN}{\mathbb{N}}
\DeclareMathOperator{\bA}{\mathbb{A}}
\DeclareMathOperator{\cH}{\mathcal{H}}

\DeclareMathOperator{\g}{\mathfrak{g}}

\DeclareMathOperator{\n}{\mathfrak{n}}
\DeclareMathOperator{\fb}{\mathfrak{b}}

\DeclareMathOperator{\cC}{\mathcal{C}}

\DeclareMathOperator{\As}{As}
\DeclareMathOperator{\cU}{\mathcal{U}}

\DeclareMathOperator{\cO}{\mathcal{O}}
\DeclareMathOperator{\cZ}{\mathcal{Z}}

\DeclareMathOperator{\cW}{\mathcal{W}}

\DeclareMathOperator{\Ind}{Ind}

\DeclareMathOperator{\Pl}{Pl}

\DeclareMathOperator{\Id}{Id}
\DeclareMathOperator{\Norm}{Norm}

\DeclareMathOperator{\cK}{\mathcal{K}}

\DeclareMathOperator{\GL}{\mathrm{GL}}
\DeclareMathOperator{\ntemp}{ntemp}
\DeclareMathOperator{\gen}{gen}
\DeclareMathOperator{\RS}{RS}
\DeclareMathOperator{\Sh}{Sh}

\numberwithin{equation}{subsection}
\newcounter{keepeqno}
\newenvironment{num}
 {\setcounter{keepeqno}{\value{equation}}%
  \begin{list}{(\theequation)}{\usecounter{equation}}%
  \setcounter{equation}{\value{keepeqno}}}
 {\end{list}}
 
\newcommand{\eq}[1][r]
   {\ar@<-3pt>@{-}[#1]
    \ar@<-1pt>@{}[#1]|<{}="gauche"
    \ar@<+0pt>@{}[#1]|-{}="milieu"
    \ar@<+1pt>@{}[#1]|>{}="droite"
    \ar@/^2pt/@{-}"gauche";"milieu"
    \ar@/_2pt/@{-}"milieu";"droite"}

\title{Archimedean theory and $\epsilon$-factors for the Asai Rankin-Selberg integrals}
\author{Rapha\"el Beuzart-Plessis \protect\footnote{The project leading to this publication has received funding from Excellence Initiative of Aix-Marseille University-A*MIDEX, a French ``Investissements d'Avenir" programme.}}

\begin{document}

\maketitle

\begin{abstract}
In this paper, we partially complete the local Rankin-Selberg theory of Asai $L$-functions and $\epsilon$-factors as introduced by Flicker and Kable. In particular, we establish the relevant local functional equation at Archimedean places and prove the equality between Rankin-Selberg's and Langlands-Shahidi's $\epsilon$-factors at every place. Our proofs work uniformly for any characteristic zero local field and use as only input the global functional equation and a globalization result for a dense subset of tempered representations that we infer from work of Finis-Lapid-M\"uller. The results of this paper are used in \cite{Beu} to establish an explicit Plancherel decomposition for $\GL_n(F)\backslash \GL_n(E)$, $E/F$ a quadratic extension of local fields, with applications to the Ichino-Ikeda and formal degree conjecture for unitary groups.
\end{abstract}

\tableofcontents

\section{Introduction}

The goal of this paper is to partially complete the local Rankin-Selberg theory of Asai $L$-functions and $\epsilon$-factors as introduced by Flicker \cite{Fli} and Kable \cite{Kab}. In particular we establish the relevant local functional equation and prove the equality between Rankin-Selberg's and Artin's $\epsilon$-factors (which are the same as Langlands-Shahidi's $\epsilon$-factors by the recent paper \cite{Shank}) in full generality (both for Archimedean and non-Archimedean fields). This will be used in the paper \cite{Beu} to establish an explicit Plancherel decomposition for $\GL_n(F)\backslash \GL_n(E)$, $E/F$ a quadratic extension of local fields, with applications to the Ichino-Ikeda and formal degree conjecture for unitary groups.

Recall that the now classical Rankin-Selberg local theory for tensor $L$-functions and $\epsilon$-factors has been developed by Jacquet-Piatetskii-Shapiro-Shalika, Jacquet-Shalika and Jacquet in \cite{JPSS}, \cite{JS2} and \cite{Jac}. The main results of those references roughly say that we can define local $L$-functions of pairs as the ``greatest common divisor" (in some loose sense) of certain families of Zeta integrals and local $\epsilon$-factors of pairs through certain functional equations satisfied by those families. Moreover, it is one of the characterizing properties of the local Langlands correspondence of \cite{HT}, \cite{Hen} and \cite{Sch} that in the $p$-adic case these so defined local $L$- and $\epsilon$-factors match the Artin $L$- and $\epsilon$-factors on the Galois side. In the Archimedean case, the Langlands correspondence is characterized by other means \cite{La} and it is a result of Jacquet and Shalika \cite{JS2}, \cite{Jac} that Artin $L$-functions of pairs can indeed be considered as the ``greatest common divisor" of the relevant family of Zeta integrals and moreover that those satisfy the correct functional equation with respect to Artin $\epsilon$-factors of pairs.

In \cite{Fli}, Flicker has introduced in the non-Archimedean case a family of Zeta integrals that ought to represent the Asai $L$-function of a given generic irreducible (smooth) representation $\pi$ of $\GL_n(E)$ where $E$ is a quadratic extension of a non-Archimedean field $F$ (and the Asai $L$-function is taken with respect to this extension). In particular, he was able to define a Rankin-Selberg type Asai $L$-function $L^{\RS}(s,\pi,\As)$ as the greatest common divisor of his family of Zeta integrals as well as a $\epsilon$-factor $\epsilon^{\RS}(s,\pi,\As,\psi')$ (where $\psi':F\to \mathbb{S}^1$ is a non-trivial character) through the existence of a functional equation satisfied by the same Zeta integrals. Similar results have been obtained independently by Kable in \cite{Kab}. To be more specific, let $\psi$ be a nontrivial additive character of $E$ which is trivial on $F$ and let $\cW(\pi,\psi)$ be the Whittaker model of $\pi$ with respect to the corresponding standard character of the standard maximal unipotent subgroup $N_n(E)$ of $\GL_n(E)$. The Zeta integrals defined by Flicker and Kable are associated to functions $W\in \cW(\pi,\psi)$ and $\phi\in C_c^\infty(F^n)$ and defined by
$$\displaystyle Z(s,W,\phi)=\int_{N_n(F)\backslash \GL_n(F)} W(h) \phi(e_nh)\lvert \det h\rvert_F^s dh$$
where $s$ is a complex parameter, $e_n=(0,\ldots,0,1)$ and $\lvert .\rvert_F$ the normalized absolute value on $F$. Flicker and Kable show that this integral converges when the real part of $s$ is large enough and that it is represented by a rational function in $q^{-s}$, $q$ being the order of the residue field of $F$. Moreover, they also prove that the vector space spanned by $\{Z(s,W,\phi)\mid W\in \cW(\pi,\psi),\phi\in C_c^\infty(F^n) \}$ is a fractional ideal for $\C[q^s,q^{-s}]$ generated by an unique element $L^{\RS}(s,\pi,\As)$ of the form $P(q^{-s})^{-1}$ where $P\in \C[T]$ is such that $P(0)=1$. The next result of \cite{Fli} and \cite{Kab} is the functional equation. For $W\in \cW(\pi,\psi)$ set $\widetilde{W}(g)=W(w_n{}^tg^{-1})$ for every $g\in \GL_n(E)$ where $w_n=\begin{pmatrix} & & 1 \\ & \iddots & \\ 1 & & \end{pmatrix}$. Notice that $\widetilde{W}\in \cW(\widetilde{\pi},\psi^{-1})$ where $\widetilde{\pi}$ stands for the contragredient of the representation $\pi$. Let $\psi':F\to \mathbb{S}^1$ be a nontrivial additive character and $\phi\mapsto \widehat{\phi}$ be the usual Fourier transform on $C_c^\infty(F^n)$ defined using the character $\psi'$ and the corresponding autodual Haar measure i.e. for every $\phi\in C_c^\infty(F^n)$ we have
$$\displaystyle \widehat{\phi}(x_1,\ldots,x_n)=\int_{F^n} \phi(y_1,\ldots,y_n)\psi'(x_1y_1+\ldots+x_ny_n)dy_1\ldots dy_n,\;\; (x_1,\ldots,x_n)\in F^n,$$
where the measure of integration is chosen such that $\widehat{\widehat{\phi}}(v)=\phi(-v)$. Then the functional equation reads as follows: there exists a unique monomial $\epsilon^{\RS}(s,\pi,\As,\psi')$ in $q^{-s}$ such that
$$\displaystyle \frac{Z(1-s,\widetilde{W},\widehat{\phi})}{L^{\RS}(1-s,\widetilde{\pi},\As)}=\omega_\pi(\tau)^{n-1}\lvert \tau\rvert_E^{\frac{n(n-1)}{2}(s-1/2)}\lambda_{E/F}(\psi')^{-\frac{n(n-1)}{2}}\epsilon^{\RS}(s,\pi,\As,\psi')\frac{Z(s,W,\phi)}{L^{\RS}(s,\pi,\As)}$$
for every $W\in \cW(\pi,\psi)$ and $\phi\in C_c^\infty(F^n)$ where $\omega_{\pi}$ stands for the central character of $\pi$, $\tau\in E$ is the unique element so that $\psi(z)=\psi'(\Tra_{E/F}(\tau z))$ for every $z\in E$ (by the assumption made on $\psi$, we have $\Tra_{E/F}(\tau)=0$) and $\lambda_{E/F}(\psi')$ is the Langlands constant of the extension $E/F$ (see Section \ref{Section Asai factors} for a reminder). The attentive reader will have noticed that the $\epsilon$-factor $\epsilon^{\RS}(s,\pi,\As,\psi')$ is normalized differently than in \cite{Fli} and \cite{Kab} (due to the appearance of the factor $\omega_\pi(\tau)^{n-1}\lvert \tau\rvert_E^{\frac{n(n-1)}{2}(s-1/2)}\lambda_{E/F}(\psi')^{-\frac{n(n-1)}{2}}$). The present normalization will be justified a posteriori by the equality between $\epsilon^{\RS}(s,\pi,\As,\psi')$ and the Artin $\epsilon$-factor $\epsilon(s,\pi,\As,\psi')$.

By the work of Anandavardhanan-Rajan \cite{AR} (for $\pi$ square-integrable) and Matringe \cite{Mat} (for general $\pi$) the $L$-function $L^{\RS}(s,\pi,\As)$ matches Shahidi's Asai $L$-function $L^{\Sh}(s,\pi,\As)$ (\cite{Sha3}, \cite{Gold}) and hence by \cite{Hen2} also the corresponding Artin $L$-function $L(s,\pi,\As)$ (defined through the local Langlands correspondence). Recently, Anandavardhanan-Kurinczuk-Matringe-S\'echerre-Stevens \cite{AKMSS} have also established the equality between $\epsilon^{\RS}(s,\pi,\As,\psi')$ with the Shahidi's $\epsilon$-factor $\epsilon^{\Sh}(s,\pi,\As,\psi')$ when $\pi$ is supercuspidal. By the recent preprint \cite{Shank} of Shankman this also gives the equality $\epsilon^{\RS}(s,\pi,\As,\psi')=\epsilon(s,\pi,\As,\psi')$ with the Artin $\epsilon$-factor when $\pi$ is supercuspidal. A similar result has been obtained when $n=2$ for any representation $\pi$ in the recent preprint \cite{CCI}. In this paper, we will complete those results in the characteristic zero case by showing that the previous equality between $\epsilon$-factors holds in general and also by working out the Archimedean theory. Moreover, we will also reprove most of the previous results in the $p$-adic case since our methods can treat uniformly the Archimedean and non-Archimedean case. Our main inputs will be the global functional equation satisfied by the corresponding global Zeta integrals (already established in \cite{Fli}, \cite{Kab}) as well as a globalization result allowing to realize a dense subset of the tempered dual of $\GL_n(E)$ as local constituents of global cuspidal automorphic representations of $\GL_n$ with a control on the ramification. We deduce this globalization result from the recent work of Finis-Lapid-M\"uller \cite{FLM} on limit multiplicities for cuspidal automorphic representations of $\GL_n$.

We now describe the main result of this paper. Let $E/F$ be a quadratic extension of local fields of characteristic zero. In the Archimedean case, by a smooth representation of $\GL_n(E)$ we will mean a smooth admissible Fr\'echet representation of moderate growth in the sense of Casselman-Wallach \cite{Cas}, \cite[Sect. 11]{Wall2} or, which is the same, an admissible SF representation in the sense of \cite{BK}. Let $\pi$ be a generic irreducible smooth representation of $\GL_n(E)$ and $\cW(\pi,\psi)$ be its Whittaker model with respect to a fixed nontrivial additive character $\psi:E\to \mathbb{S}^1$ which we again take to be trivial on $F$. To $W\in \cW(\pi,\psi)$ we associate $\widetilde{W}\in \cW(\widetilde{\pi},\psi^{-1})$ as before. Let $\cS(F^n)$ be $C_c^\infty(F^n)$ in the $p$-adic case and the usual Schwartz space in the Archimedean case. We let $\phi\mapsto \widehat{\phi}$ be the usual Fourier transform on $\cS(F^n)$ defined using a nontrivial additive character $\psi':F\to \mathbb{S}^1$ as before. Let $\tau\in E$ be the unique element so that $\psi(z)=\psi'(\Tra_{E/F}(\tau z))$ for every $z\in E$. Then, for $W\in \cW(\pi,\psi)$ and $\phi\in \cS(F^n)$ we define as above, whenever convergent, a Zeta integral $Z(s,W,\phi)$. The main result of this paper can now be stated as follows (see Theorems \ref{theo 1 inert} and \ref{theo 2 inert}):

\begin{theointro}\label{theo 1 intro}
Let $W\in \cW(\pi,\psi)$ and $\phi\in \cS(F^n)$. Then:
\begin{enumerate}[(i)]
\item The integral defining $Z(s,W,\phi)$ is convergent when the real part of $s$ is sufficiently large and moreover it extends to a meromorphic function on $\C$.
\item We have the functional equation
$$\displaystyle \frac{Z(1-s,\widetilde{W},\widehat{\phi})}{L(1-s,\widetilde{\pi},\As)}=\omega_\pi(\tau)^{n-1}\lvert \tau\rvert_E^{\frac{n(n-1)}{2}(s-1/2)}\lambda_{E/F}(\psi')^{-\frac{n(n-1)}{2}}\epsilon(s,\pi,\As,\psi')\frac{Z(s,W,\phi)}{L(s,\pi,\As)}$$
where $L(s,\widetilde{\pi},\As)$, $L(s,\pi,\As)$ and $\epsilon(s,\pi,\As,\psi')$ stand for the Asai $L$- and $\epsilon$-factors of Artin.
\item The function $\displaystyle s\mapsto \frac{Z(s,W,\phi)}{L(s,\pi,\As)}$ is holomorphic. Moreover, if $\pi$ is {\em nearly tempered} (see Section \ref{Part II}), for every $s_0\in \C$ we can choose $W\in \cW(\pi,\psi)$ and $\phi\in \cS(F^n)$ such that this function does not vanish at $s_0$.
\end{enumerate}
\end{theointro}

We recall here that, in the $p$-adic case, the theorem above has already been proved by Flicker \cite[Appendix]{Fli} and Kable \cite[Theorem 3]{Kab} without the restriction on $\pi$ in the last statement of (iii) but with different definitions of the $L$ and $\epsilon$-factors that we denote here by $L^{\RS}(s,\pi,\As)$ and $\epsilon^{\RS}(s,\pi,\As,\psi')$. Then, Anandavardhanan-Rajan \cite{AR} (for discrete series) and Matringe \cite{Mat} (extending the result in general) have shown that $L^{\RS}(s,\pi,\As)=L(s,\pi,\As)$. Finally, the equality $\epsilon^{\RS}(s,\pi,\As,\psi')=\epsilon(s,\pi,\As,\psi')$ was already known when $\pi$ is supercuspidal (\cite{AKMSS}) or $n=2$ (\cite{CCI}). Therefore, in the $p$-adic case, the only new aspect of the above theorem is the equality of $\epsilon$-factors $\epsilon^{\RS}(s,\pi,\As,\psi')=\epsilon(s,\pi,\As,\psi')$ in general. On the other hand, the Archimedean case seems to have remained completely untouched.

One comment is in order: we expect $L(s,\pi,\As)$ to be the ``greatest common divisor'', in a suitable sense, of the family of Zeta integrals $Z(s,W,\phi)$ and therefore the second part of (iii) should be true in general. This is known in the $p$-adic case by the aforementioned result of Matringe but remains to be seen in the Archimedean case for irreducible generic representations $\pi$ which are not nearly tempered in the sense of Section \ref{Part II}. Actually, we even expect something stronger: there should exist finite families $W_i\in \cW(\pi,\psi)$ and $\phi_i\in \cS(F^n)$ such that
$$\displaystyle L(s,\pi,\As)=\sum_i Z(s,W_i,\phi_i).$$
Again in the $p$-adic case, this is known by the very definition of $L^{\RS}(s,\pi,\As)$ as a GCD and the result of Matringe. However, in the Archimedean case this result seems harder to establish and in any case unreachable by the method developed in this paper. We refer the reader to \cite[Theorem 2.7]{Jac} for a proof of this property in the case of the usual Rankin-Selberg integrals.

We now briefly describe the content of each section of this paper. In Part \ref{Part preliminaries} we gather general results which are not specific to $\GL_n$. In Section \ref{Section groups}, we set up most of our notation for a general connected reductive group. Section \ref{Section TVS} contains a reminder on properties  of topological vector spaces and functions valued in them that we shall use repeatedly in the paper. In Section \ref{Section representations}, we set up some notation and conventions related to representations of local reductive groups. Section \ref{Section matrix coefficients} contains a probably well-known result giving uniform bounds for matrix coefficients of a family of parabolically induced representations. In view of the lack of a proper reference, we provide a proof. In Section \ref{Section HCS spaces}, we introduce various ``extended" Harish-Chandra Schwartz spaces associated to a generic character on a quasi-split local reductive group. These spaces will be the natural receptacle for the Whittaker models of a family of induced representations. Section \ref{Section bounds Whittaker} is devoted to establishing locally uniform bounds for analytic family of Whittaker functions living in the Whittaker models of a family of parabolically induced representations. The proof proceeds by a reduction to the case of usual matrix coefficients by ``smoothing" the Whittaker functional and is far more technical in the Archimedean case. Similar (point-wise) bounds have been obtained by Wallach \cite[Theorem 15.2.5]{Wall2}. Using the result of Section \ref{Section bounds Whittaker} and the theory of the Jacquet's functional, we construct in Section \ref{Section good sections Whittaker} certain ``good" sections for Whittaker models of a family of parabolically induced representations. The existence of such families will be a crucial ingredient to extend Theorem \ref{theo 1 intro} from a dense subset of the tempered dual to the set of all generic representations. In Section \ref{Section automatic holomorphic cont}, we establish a result on the automatic holomorphic continuation of certain functions in several complex variables which are of ``finite order'' in vertical strips in one of the variable locally uniformly in the remaining ones and satisfying a certain functional equation. This will later be applied to show the meromorphic continuation of the local Asai Rankin-Selberg integrals and control their poles in terms of Asai $L$-factors. The theory of local Asai Zeta integrals and their functional equations is the object of Part \ref{Part II}. Section \ref{Section Asai factors} is a reminder on basic properties of local Asai $L$- and $\epsilon$-factors of Artin type. In Section \ref{Section defn and convergence}, we introduce the relevant Zeta integrals and establish basic convergence results on them. In Section \ref{Section split case}, we recall the meromorphic continuation and functional equation of these integrals in the split case (i.e. when $E=F\times F$) which is due to Jacquet, Piatetski-Shapiro and Shalika \cite{JPSS}, \cite{JS2}, \cite{Jac}. In Section \ref{Section inert case}, we state the main theorems of this paper pertaining to the same meromorphic continuation and functional equation but in the inert case (i.e. when $E/F$ is a field extension). Section \ref{Section unr computation} is devoted to the unramified computation of these Zeta integrals which is already in the literature in all but one case (i.e. when $E/F$ is a ramified quadratic extension but $\pi$, $\psi'$ and $\psi$ are unramified). In Section \ref{Global functional equation}, we recall the definition of the global Zeta integrals and their functional equation. Section \ref{Section globalization} contains the aforementioned globalization result due to Finis-Lapid-M\"uller. The deduction from \cite{FLM} is carefully explained. We also sketch how a similar result (in a slightly weaker form) can be proved for general reductive groups using part of {\it loc. cit.}. Finally, Sections \ref{appendix proof of the theorem} and \ref{Section end of proof of main theorems} contain the proof of the main results.

\subsection{General notation}

In this paper $F$ will always be a local field of characteristic zero. We will denote by $\lvert .\rvert_F$ the normalized absolute value of $F$. In the non-Archimedean case, we let $\cO_F$ be the ring of integers of $F$ and $q=q_F$ be the cardinality of the residue field of $F$.

For two positive functions $f_1$, $f_2$ on a set $X$ a sentence like
\begin{center}
$f_1(x)\ll f_2(x)$ for all $x\in X$
\end{center}
means that there exists a constant $C>0$ such that $f_1(x)\leqslant Cf_2(x)$ for every $x\in X$. When we want to emphasize that the implicit constant depends on some auxiliary parameters $a,b,c...$ we will write ``$f_1(x)\ll_{a,b,c...} f_2(x)$ for all $x\in X$".

For every complex number $z\in \C$, we write $\Re(z)$ and $\Im(z)$ for the real and imaginary parts of $z$ respectively.

\section{Preliminaries}\label{Part preliminaries}

\subsection{Groups}\label{Section groups}

Let $G$ be a connected reductive group over $F$. We denote by $A_G$ the maximal split torus in the center of $G$ and by $X^*(G)$, $X^*(A_G)$ the groups of algebraic characters of $G$ and $A_G$ respectively. Set
$$\displaystyle \cA_G^*=X^*(G)\otimes \bR=X^*(A_G)\otimes \bR,\;\;\; \cA_{G,\C}^*=X^*(G)\otimes \C=X^*(A_G)\otimes \C$$
and
$$\displaystyle \cA_G=\Hom(X^*(G),\bR),\;\;\; \cA_{G,\C}=\Hom(X^*(G),\C)$$
for their duals. We denote by $\langle .,.\rangle$ the natural pairing between $\cA_G^*$ and $\cA_G$ (resp. between $\cA_{G,\C}^*$ and $\cA_{G,\C}$). Let $H_G:G(F)\to \cA_G$ be the homomorphism characterized by $\langle \chi,H_G(g)\rangle=\log\lvert \chi(g)\rvert_F$ for every $g\in G(F)$ and $\chi\in X^*(G)$. For $\lambda\in \cA_{G,\C}^*$ we denote by $g\mapsto g^\lambda$ the unramified character of $G(F)$ given by $g^\lambda=e^{\langle \lambda,H_G(g)\rangle}$. Moreover, we let $\Re(\lambda)\in \cA_G^*$ be the {\em real part} of $\lambda$ (i.e. its projection to $\cA_G^*$ relative to the decomposition $\cA_{G,\C}^*=\cA_G^*\oplus i\cA_G^*$). Let $P=MN$ be a parabolic subgroup of $G$. Then the restriction map $X^*(A_M)\to X^*(A_G)$ induces surjections $\cA_M^*\to \cA_G^*$, $\cA_{M,\C}^*\to \cA_{G,\C}^*$ whose kernels will be denoted by $(\cA_M^G)^*$ and $(\cA_{M,\C}^G)^*$ respectively.

We fix a minimal parabolic subgroup $P_0$ of $G$ with a Levi decomposition $P_0=M_0N_0$ and we set $A_0=A_{M_0}$, $\cA_0^*=\cA_{M_0}^*$, $(\cA^G_0)^*=(\cA^G_{M_0})^*$, $\cA_0=\cA_{M_0}$, $H_0=H_{M_0}$. We also choose a maximal compact subgroup $K$ of $G(F)$ which is special in the $p$-adic case and in good position relative to $M_0$. We endow $K$ with its unique Haar measure of total mass one. For every parabolic subgroup $P$ of $G$ we have the Iwasawa decomposition $G(F)=P(F)K$. As usual, by a {\em standard} parabolic subgroup we mean a parabolic subgroup of $G$ containing $P_0$. If $P$ is such a standard parabolic subgroup we will always write $P=MN$ for its unique Levi decomposition with $M_0\subseteq M$ and $\overline{P}=M\overline{N}$ for the opposite parabolic subgroup. The restriction map $X^*(M)\to X^*(M_0)$ then induces an embedding $\cA_M^*\hookrightarrow \cA_0^*$ through which we will always consider $\cA_M^*$ as a subspace of $\cA_0^*$. We shall also denote by $\delta_P$ the modular character of $P(F)$. By the choice of $K$, for every standard parabolic subgroup $P=MN$ we have $K_P=K_MK_N$ where $K_P=K\cap P(F)$, $K_M=K\cap M(F)$ and $K_N=K\cap N(F)$.

Let $\Delta\subseteq X^*(A_0)$ be the set of simple roots of $A_0$ in $N_0$ and $\Delta^\vee\subseteq \cA_0$ the corresponding subset of simple coroots. We set
$$\displaystyle \overline{\cA_0^+}=\{ X\in \cA_0\mid \langle X, \alpha\rangle \leqslant 0 \; \forall \alpha\in \Delta\}$$
$$\displaystyle \overline{(\cA_0^*)^+}=\{ \lambda\in \cA_0^*\mid \langle \lambda,\alpha^\vee\rangle\leqslant 0 \; \forall \alpha\in \Delta\}$$
for the {\em closed negative Weyl chambers} in $\cA_0$ and $\cA_0^*$ respectively. We also let
$$\displaystyle M_0^+=H_0^{-1}(\overline{\cA_0^+})=\{ m_0\in M_0(F)\mid m_0^\alpha\leqslant 1, \; \forall \alpha\in \Delta\}.$$
Let $W^G=\Norm_{G(F)}(M_0)/M_0(F)$ be the Weyl group of $M_0$. Then, $W^G$ acts on $\cA_0^*$ and $\overline{(\cA_0^*)^+}$ is a fundamental domain for this action. For every $\lambda\in \cA_0^*$ we denote by $\lvert \lambda \rvert$ the unique element in the intersection $W^G\lambda\cap \overline{(\cA_0^*)^+}$. We equip $\cA_0^*$ with the strict partial order $\prec$ defined by
\begin{center}
$\lambda\prec \mu$ if and only if $\mu-\lambda=\sum_{\alpha\in \Delta} x_\alpha \alpha$ where $x_\alpha>0$ for every $\alpha\in \Delta$.
\end{center}

We fix an algebraic group embedding $\iota:G/A_G\hookrightarrow GL_N$ for some $N\geqslant 1$ and for every $g\in G(F)$ we set
$$\displaystyle \overline{\sigma}(g)=\sup\left(\{1\}\cup \{\log \lvert \iota(g)_{i,j}\rvert_F\mid 1\leqslant i,j\leqslant N \}\right)$$
where the $\iota(g)_{i,j}$'s denote the entries of the matrix $\iota(g)$.

We denote by $\mathfrak{g}$ the (algebraic) Lie algebra of $G$. Similar notations will be used for other algebraic groups (i.e. we will denote the Lie algebra of an algebraic group by the corresponding gothic letter). In the Archimedean case, we will also write $\mathfrak{k}$ for the (real) Lie algebra of $K$ and $\cU(\mathfrak{k})$, $\cU(\mathfrak{g})$ for the enveloping algebras of the complexifications of $\mathfrak{k}$ and $\mathfrak{g}(F)$ (considered as a real Lie algebra) respectively. We identify every element of $\cU(\g)$ with the distribution supported at $1$ that it defines.

We will assume that all the locally compact topological groups that we encounter have been equipped with Haar measures (bi-invariant Haar measures as we will always integrate, with one exception, over unimodular groups). The precise choices of these Haar measures will always be irrelevant. We denote by $\ast$ the convolution product on a locally compact topological group $H$. In the Archimedean case, this convolution product extends to distributions of compact support including elements of $\cU(\g)$ and continuous compactly supported functions on closed subgroups (seen as distributions through the choice of a Haar measure on that subgroup).

\subsection{Topological vector spaces}\label{Section TVS}

In this paper, by a {\em topological vector space} (TVS) we always mean a Hausdorff locally convex topological vector space over $\C$. If $F$ is a TVS, we shall denote by $F'$ its continuous dual. The TVS to be considered in this paper will all be LF spaces that is countable direct limit of Fr\'echet spaces. We will even only encounter {\em strict} LF spaces i.e. TVS $F$ that can be written as the direct limit of a sequence $(F_n)_n$ of Fr\'echet spaces where the transition maps $F_n\to F_{n+1}$ are closed embeddings. Strict LF spaces are complete. Moreover, since LF spaces are barreled \cite[Corollary 3 of Proposition 33.2]{Tr}, they satisfy the uniform boundedness principle (aka Banach-Steinhaus theorem).

We refer the reader to \cite{Bour} for the notions of smooth and holomorphic functions valued in topological vector spaces (TVS) that will be used thoroughly in this paper. Actually, we will only consider smooth and holomorphic functions valued in Fr\'echet or strict LF spaces for which the following criterion can be applied (\cite[3.3.1 (v)]{Bour}):
\begin{num}
\item\label{eq 1 TVS} Let $F$ be a quasi-complete TVS and $M$ a complex analytic manifold. Then, a function $\varphi:M\to F$ is holomorphic if and only if it is continuous and for some total subspace $H\subset F'$ the scalar-valued functions $m\in M\mapsto \langle \varphi(m),\lambda\rangle$ are holomorphic for every $\lambda\in H$.
\end{num}

Also when we say that a function on a totally disconnected locally compact topological space (e.g. $G(F)$ in the $p$-adic case) is {\em smooth} we always mean that it is locally constant.

\subsection{Representations}\label{Section representations}

By a {\em representation} of $G(F)$ we will always mean a smooth representation of finite length with complex coefficients. Here {\em smooth} has the usual meaning in the $p$-adic case (i.e. every vector has an open stabilizer) whereas in the Archimedean case it means a smooth admissible Fr\'echet representation of moderate growth in the sense of Casselman-Wallach \cite{Cas}, \cite[Sect. 11]{Wall2} or, which is the same, an admissible SF representation in the sense of \cite{BK}. We shall always abuse notation and denote by the same letter a representation and the space on which it acts. In the Archimedean case this space always comes with a topology (it is a Fr\'echet space) whereas in the $p$-adic case it will sometimes be convenient, in order to make uniform statements, to equip this space with its finest locally convex topology (it then becomes a strict LF space). For $\pi$ a representation of $G(F)$ and $\lambda\in \cA_{G,\C}^*$, we denote by $\pi_\lambda$ the twist of $\pi$ by the character $g\in G(F)\mapsto g^\lambda$. We let $\Irr(G)$, $\Temp(G)$ and $\Pi_2(G)$ be the sets of isomorphism classes of all irreducible representations, irreducible tempered representations and irreducible square-integrable representations of $G(F)$ respectively. We denote by $\widetilde{\pi}$ the contragredient representation of $\pi$ (aka smooth dual) and by $\langle .,.\rangle$ the natural pairing between $\pi$ and $\widetilde{\pi}$. In the Archimedean case, $\widetilde{\pi}$ can be defined as the Casselman-Wallach globalization of the contragredient of the Harish-Chandra module underlying $\pi$. An alternative description of $\widetilde{\pi}$, which is more suitable in practice, is as the space of linear forms on $\pi$ which are continuous with respect {\em to any} $G(F)$-continuous norm on $\pi$ together with the natural $G(F)$-action on it (a norm on $\pi$ is said to be {\em $G(F)$-continuous} if the action of $G(F)$ on $\pi$ is continuous for this norm). In the Archimedean case, we denote by $\pi'$ the topological dual of $\pi$ (i.e. the space of all continuous linear forms on $\pi$). There is a natural action $g\mapsto \pi'(g)$ of $G(F)$ on $\pi'$ and if we equip $\pi'$ with the topology of uniform convergence on compact subsets of $\pi$ then this action of $G(F)$ on $\pi'$ is continuous and this allows one to define $\pi'(\varphi)$ for every function $\varphi\in C_c(G(F))$ by integration.

Let $P=MU$ be a parabolic subgroup of $G$. For $\sigma$ a representation of $M(F)$, we denote by $i_P^G(\sigma)$ the smooth normalized parabolic induction of $\sigma$ from $P(F)$ to $G(F)$. The space of $i_P^G(\sigma)$ consists of smooth functions $e:G(F)\to \sigma$ satisfying $e(mug)=\delta_P(m)^{1/2}\sigma(m)e(g)$ for every $(m,u,g)\in M(F)\times U(F)\times G(F)$ (with its natural structure of Fr\'echet space in the Archimedean case) and the group $G(F)$ acts by right translation. Assume that $P$ is standard. Then, for every $\lambda\in \cA_{M,\C}^*$ restriction to $K$ induces a topological isomorphism between $\pi_\lambda=i_P^G(\sigma_\lambda)$ and $\pi_K=i_{K_P}^K(\sigma_{\mid K_M})$ where $K_P=K\cap P(F)$, $K_M=K\cap M(F)$ and $i_{K_P}^K(\sigma_{\mid K_M})$ denotes the spaces of smooth functions $e:K\to \sigma$ satisfying $e(muk)=\sigma(m)e(k)$ for every $(m,u,k)\in K_M\times (K\cap U(F))\times K$. These identifications allows to define a notion of {\em holomorphic sections} $\lambda\in \cA_{M,\C}^*\mapsto e_\lambda\in \pi_{\lambda}$: we call such a map holomorphic if its composition with the isomorphism $\pi_{\lambda}\simeq \pi_K$ is holomorphic. It turns out that this notion is actually independent of the choice of $K$. We also say that an assignment $\lambda\in \cA_{M,\C}^*\mapsto T_\lambda\in \pi_{\lambda}'$ is {\em holomorphic} if for every holomorphic section $\lambda\mapsto e_\lambda\in \pi_{\lambda}$ the map $\lambda\in \cA_{M,\C}^*\mapsto T_\lambda(e_\lambda)$ is holomorphic. It is equivalent to ask that, identifying $T_\lambda$ with a linear form on $\pi_K$ for every $\lambda\in \cA_{M,\C}^*$, for every $e\in \pi_K$ the map $\lambda\in \cA_{M,\C}^*\mapsto T_\lambda(e)$ be holomorphic: Indeed as $\pi_K$ is Fr\'echet hence barreled, if this is so by the Banach-Steinhaus theorem for every compact $\cK\subseteq \cA_{M,\C}^*$ the subset $\{T_\lambda\mid \lambda\in \cK \}$ of $\pi_K'$ is equicontinuous hence $\cA_{M,\C}^*\times \pi_K\to \C$, $(\lambda,e)\mapsto T_\lambda(e)$ is continuous and then we can apply \cite[Lemma 1]{Jac2}.

\begin{lem}\label{lem 1 representations}
Assume that $F$ is Archimedean. Let $\pi\in \Irr(G)$ and let $C_K$ be a Casimir element for $K$ (i.e. the element of $\cU(\mathfrak{k})$ associated to a negative definite $K$-invariant real symmetric bilinear form on the Lie algebra of $K$). Then,
\begin{enumerate}[(i)]
\item $\pi(1+C_K)$ is a topological isomorphism of $\pi$ onto itself.
\end{enumerate}
Let $p_0$ be a $G(F)$-continuous Hilbert norm (i.e. associated to a scalar product) on $\pi$ which is $K$-invariant and set $p_\ell=p_0\circ \pi(1+C_K)^\ell$ for every integer $\ell\in \bZ$. Then,
\begin{enumerate}[(i)]
\setcounter{enumi}{1}
\item The map $\ell\mapsto p_\ell$ is increasing i.e. $p_\ell(v)\leqslant p_{\ell+1}(v)$ for every $\ell\in \bZ$ and $v\in \pi$.
\item The family of norms $(p_\ell)_{\ell \geqslant 0}$ generates the topology on $\pi$.
\item Let $\pi^{(\ell)}$ be the completion of $\pi$ for $p_\ell$. Then, the natural pairing $\pi\times \widetilde{\pi}\to \C$ extends to a continuous bilinear form $\pi^{(\ell)}\times \widetilde{\pi}\to \C$ giving an embedding $\pi^{(\ell)}\subset (\widetilde{\pi})'$ for every $\ell$ and any equicontinuous subset of $(\widetilde{\pi})'$ is contained and bounded in $\pi^{(\ell)}$ for some $\ell$.
\end{enumerate}
\end{lem}

\noindent\ul{Proof}: For each $\gamma\in \widehat{K}$ the element $1+C_K$ acts on $\gamma$ by a scalar $N(\gamma)$ which is greater or equal to $1$ (Indeed, $C_K$ being in $\cU(\mathfrak{k})^K$ the corresponding operator on $\gamma$ is scalar and as $C_K=-\sum_i X_i^2$ for a certain basis $(X_i)$ of $\mathfrak{k}$ it is also positive hermitian with respect to any $K$-invariant scalar product on the space of $\gamma$). Let $p_0$ be a $G(F)$-continuous Hilbert norm on $\pi$ which is $K$-invariant (such norm exists as the Harish-Chandra module underlying $\pi$ admits at least one Hilbert globalization). Set $p_\ell=p_0\circ \pi(1+C_K)^\ell$ for every integer $\ell\geqslant 0$. Then, by \cite[Proposition 3.9]{BK} the family of norms $(p_\ell)_{\ell \geqslant 0}$ generates the topology on $\pi$. This gives (iii). For every vector $v\in \pi$ write $\displaystyle v=\sum_{\gamma\in \widehat{K}}v_\gamma$ for its decomposition in $K$-isotypic components (the corresponding series is absolutely convergent in $\pi$). Then, by the $K$-invariance of $p_0$ for every $\ell\geqslant 0$ and $v\in \pi$ we have
\[\begin{aligned}
\displaystyle p_{\ell+1}(v) & =\left(\sum_{\gamma\in \widehat{K}}p_0(\pi(1+C_K)^{\ell+1} v_\gamma)^2 \right)^{1/2}=\left(\sum_{\gamma\in \widehat{K}}N(\gamma)^{2\ell+2}p_0(v_\gamma) \right)^{1/2} \\
 & \geqslant \left(\sum_{\gamma\in \widehat{K}}N(\gamma)^{2\ell}p_0(v_\gamma) \right)^{1/2}=p_\ell(v).
\end{aligned}\]
As $p_\ell(\pi(1+C_K)v)=p_{\ell+1}(v)$ this shows that $\pi(1+C_K)$ realizes a topological isomorphism between $\pi$ and one of its closed subspace. Since the subspace of $K$-finite vectors of $\pi$ is obviously contained in the image of $\pi(1+C_K)$ and is dense in $\pi$ this shows (i) and therefore the definition of $p_\ell$ now also makes sense for $\ell\leqslant 0$. Moreover, the above inequality still holds for every $\ell\in \bZ$ (with the same proof). This gives (ii). Finally, let $\widetilde{p}_0$ be the norm on $\widetilde{\pi}$ dual to $p_0$ i.e.
$$\displaystyle \widetilde{p}_0(v^\vee)=\sup_{\substack{v\in \pi \\ p_0(v)=1}} \lvert \langle v,v^\vee\rangle \rvert,\;\;\; v^\vee\in \widetilde{\pi}.$$
Set $\widetilde{p}_\ell=\widetilde{p}_0\circ \widetilde{\pi}(1+C_K)^\ell$ for every $\ell\in \bZ$. Then, we easily check that $p_{-\ell}$ is the norm dual to $\widetilde{p}_\ell$ for every $\ell\in \bZ$. This gives (iv) since the family of norms $(\widetilde{p}_\ell)_{\ell\geqslant 0}$ generates the topology on $\widetilde{\pi}$. $\blacksquare$

\subsection{Uniform bounds for matrix coefficients}\label{Section matrix coefficients}

For every $\lambda\in \cA_0^*$, set $\pi^0_\lambda=i_{P_0}^G(\lambda)$ where we have identified $\lambda$ with the character it defines on $M_0(F)$. We identify the space of $\pi^0_\lambda$ with $\pi^0_K:=i_{K_{P_0}}^K(1)$ by restriction to $K$ and we equip it with the $K$-invariant pairing given by
$$\displaystyle \langle e,e'\rangle_0=\int_K e(k)e'(k)dk,\;\;\; e,e'\in \pi^0_K.$$
By the previous identifications, $\langle .,.\rangle_0$ induces a $G(F)$-invariant continuous pairing between $\pi^0_\lambda$ and $\pi^0_{-\lambda}$ for every $\lambda\in \cA_0^*$. Let $e_0\in \pi^0_K$ be the vector defined by $e_0(k)=1$ for every $k\in K$. Then, we let
$$\displaystyle \Xi^G_\lambda(g)=\langle \pi^0_\lambda(g)e_0,e_0\rangle_0,\;\;\; \lambda\in \cA_0^*, g\in G(F).$$
When $\lambda=0$, we simply set $\Xi^G=\Xi^G_0$ (it is the usual Harish-Chandra's Xi function see \cite[Sect. II.8.5]{Var} and \cite[Sect. II.1]{Wald1}). We summarize the basic properties of the functions $\Xi^G_\lambda$ in the next proposition. In the Archimedean case, most of this is contained in \cite[Sect. 3.6]{Wall}.

\begin{prop}\label{prop 1 uniform bounds matrix coefficients}
\begin{enumerate}[(i)]
\item $\Xi^G_{w\lambda}=\Xi^G_\lambda$ for every $\lambda\in \cA_0^*$ and every $w\in W^G$.

\item For every $\lambda\in \cA_0^*$, $\mu\in \cA_G^*$ and $g\in G(F)$, we have $\Xi^G_{\lambda+\mu}(g)=\Xi^G_\lambda(g)g^{\mu}$.

\item Let $P=MU$ be a parabolic subgroup containing $P_0$ (with the Levi component chosen so that $M\supset M_0$). Choose for every $g\in G(F)$ a decomposition $g=m_P(g)u_P(g)k_P(g)$ with $(m_P(g),u_P(g),k_P(g))\in M(F)\times U(F)\times K$. Then, we have
$$\displaystyle \Xi^G_\lambda(g)=\int_K \delta_P(m_P(kg))^{1/2} \Xi^M_\lambda(m_P(kg))dk$$
for every $\lambda\in \cA_0^*$ and $g\in G(F)$.

\item There exists $d>0$ such that
$$\displaystyle \delta_0(m_0)^{1/2}m_0^{\lvert \lambda\rvert}\ll \Xi^G_\lambda(m_0)\ll \delta_0(m_0)^{1/2}m_0^{\lvert \lambda\rvert} \overline{\sigma}(m_0)^d$$
for every $\lambda\in \cA_0^*$ and $m_0\in M_0^+$ where we have set $\delta_0=\delta_{P_0}$.
\end{enumerate}
\end{prop}

\noindent\ul{Proof}:
\begin{enumerate}[(i)]
\item This is \cite[Proposition 3.6.2]{Wall} in the Archimedean case and \cite[Proposition 4.1]{Cas2} in the non-Archimedan case.
\item is straightforward.
\item follows readily from the isomorphism of ``induction by stages" $\pi_\lambda^0\simeq i_P^G(i_{P_0\cap M}^M(\lambda))$ (see \cite[Lemme II.1.6]{Wald1} for the case $\lambda=0$).
\item First we prove the lower bound. By (i), we may assume $\lambda=\lvert \lambda\rvert$. Let $\overline{P}_0=M_0\overline{N}_0$ be the parabolic subgroup opposite to $P_0$ and $C\subseteq N_0(F)M_0(F)\overline{N}_0(F)$ a compact neighborhood of $1$. We can find a compact subset $C_0\subseteq M_0(F)$ such that $m_{P_0}(km_0)\in m_0C_0$ for every $k\in C$ and $m_0\in M_0^+$. By (iii) and (ii), it follows that
\[\begin{aligned}
\displaystyle \Xi^G_\lambda(m_0)=\int_K \delta_0(m_{P_0}(km_0))^{1/2} m_{P_0}(km_0)^\lambda dk & \geqslant \int_{K\cap C} \delta_0(m_{P_0}(km_0))^{1/2} m_{P_0}(km_0)^\lambda dk \\
 & \gg \delta_0(m_0)^{1/2} m_0^\lambda=\delta_0(m_0)^{1/2} m_0^{\lvert \lambda\rvert }
\end{aligned}\]
for all $m_0\in M_0^+$. We now prove the upper bound. By \cite{Kos} and \cite[Proposition 4.4.4]{BT} (see also \cite{Sil}), for every $k\in K$ and $m_0\in M_0(F)$, $H_{0}(m_{P_0}(km_0))$ belongs to the convex hull of $\{wH_0(m_0)\mid w\in W^G \}$. It follows that $m_{P_0}(km_0)^\lambda\leqslant m_0^{\lvert \lambda\rvert}$ for every $k\in K$ and $m_0\in M_0^+$. Therefore, using (ii) and (iii), we obtain
\[\begin{aligned}
\displaystyle \Xi^G_\lambda(m_0)=\int_K \delta_0(m_{P_0}(km_0))^{1/2} m_{P_0}(km_0)^\lambda dk\leqslant m_0^{\lvert \lambda\rvert}\int_K \delta_0(m_{P_0}(km_0))^{1/2} dk=m_0^{\lvert \lambda\rvert} \Xi^G(m_0)
\end{aligned}\]
for all $m_0\in M_0^+$. On the other hand, by \cite[Theorem 30 p.339]{Var} and \cite[Lemme II.1.1]{Wald1} there exists $d>0$ such that $\Xi^G(m_0)\ll \delta_0(m_0)^{1/2}\overline{\sigma}(m_0)^d$ for all $m_0\in M_0^+$. The upper bound follows. $\blacksquare$
\end{enumerate}

Let $P=MU$ be a standard parabolic subgroup and $\sigma\in \Temp(M)$. For every $\lambda\in \cA_{M,\C}^*$, we set $\pi_\lambda=i_P^G(\sigma_{\lambda})$ and identify the underlying space with $\pi_K:=i_{K_P}^K(\sigma_{\mid K_M})$ where we have set $K_P=K\cap P(F)$ and $K_M=K\cap M(F)$ as before. Similarly, we let $\widetilde{\pi}_\lambda=i_P^G(\widetilde{\sigma}_{\lambda})$ and identify its underlying space with $\widetilde{\pi}_K:=i_{K_P}^K(\widetilde{\sigma}_{\mid K_M})$ for every $\lambda\in \cA_{M,\C}^*$. We define a bilinear pairing between $\pi_K$ and $\widetilde{\pi}_K$ by
$$\displaystyle \langle e,e^\vee\rangle=\int_K \langle e(k),e^\vee(k)\rangle_\sigma dk,\;\;\; e\in \pi_K, e^\vee \in \widetilde{\pi}_K$$
where $\langle .,.\rangle_\sigma$ denotes the natural pairing between $\sigma$ and $\widetilde{\sigma}$. By the previous identifications, $\langle .,.\rangle$ induces a continuous $G(F)$-invariant pairing between $\pi_{\lambda}$ and $\widetilde{\pi}_{-\lambda}$ which identifies the latter with the smooth contragredient of $\pi_{\lambda}$ for every $\lambda\in \cA_{M,\C}^*$. The following proposition gives uniform bounds for the matrix coefficients of the $\pi_\lambda$. It is most probably well-known but in view of the lack of a proper reference, we include a proof (see however \cite[Proposition 7.14]{Knapp} for the case of $K$-finite coefficients in the Archimedean case).

\begin{prop}\label{prop 2 uniform bounds matrix coefficients}
There exist continuous semi-norms $p$ and $\widetilde{p}$ on $\pi_K$ and $\widetilde{\pi}_K$ such that
$$\displaystyle \lvert \langle\pi_\lambda(g)e,e^\vee\rangle\rvert \leqslant \Xi^G_{\Re(\lambda)}(g)p(e)\widetilde{p}(e^\vee)$$
for every $\lambda\in \cA_{M,\C}^*$, $g\in G(F)$ and $(e,e^\vee)\in \pi_K\times \widetilde{\pi}_K$.
\end{prop}

\noindent\ul{Proof}: We have
\[\begin{aligned}
\displaystyle \langle\pi_\lambda(g)e,e^\vee\rangle=\int_K \delta_P(m_P(kg))^{1/2} m_P(kg)^\lambda \langle\sigma(m_P(kg))e(k_P(kg)),e^\vee(k)\rangle_\sigma dk
\end{aligned}\]
and therefore
\[\begin{aligned}
\displaystyle \lvert \langle\pi_\lambda(g)e,e^\vee\rangle\rvert \leqslant \int_K \delta_P(m_P(kg))^{1/2} m_P(kg)^{\Re(\lambda)} \lvert \langle \sigma(m_P(kg))e(k_P(kg)),e^\vee(k)\rangle_\sigma\rvert dk
\end{aligned}\]
for every $\lambda\in \cA_{M,\C}^*$, $g\in G(F)$ and $(e,e^\vee)\in \pi_K\times \widetilde{\pi}_K$. By \cite[Theorem 2]{CHH} and \cite{Sun}, there exists continuous semi-norms $q$, $\widetilde{q}$ on $\sigma$ and $\widetilde{\sigma}$ such that $\lvert \langle \sigma(m)v,v^\vee\rangle\rvert \leqslant \Xi^M(m)q(v)\widetilde{q}(v^\vee)$ for every $m\in M(F)$ and $(v,v^\vee)\in \sigma\times \widetilde{\sigma}$. (Here we emphasize that {\em any} semi-norm on $\sigma$ or $\widetilde{\sigma}$ is continuous in the $p$-adic case). Consequently,
\[\begin{aligned}
\displaystyle \lvert \langle \pi_\lambda(g)e,e^\vee \rangle\rvert \leqslant p(e)\widetilde{p}(e^\vee)\int_K \delta_P(m_P(kg))^{1/2} m_P(kg)^{\Re(\lambda)}  \Xi^M(m_P(kg)) dk 
\end{aligned}\]
for every $\lambda\in \cA_{M,\C}^*$, $g\in G(F)$, $(e,e^\vee)\in \pi_K\times \widetilde{\pi}_K$ where we have set $p(e)=\sup_{k\in K} q(e(k))$ and $\widetilde{p}(e^\vee)=\sup_{k\in K} \widetilde{q}(e^\vee(k))$. By Proposition \ref{prop 1 uniform bounds matrix coefficients}, the above integral is equal to $\Xi^G_{\Re(\lambda)}(g)$ whereas $p$ and $\widetilde{p}$ clearly define continuous semi-norms on $\pi_K$ and $\widetilde{\pi}_K$ respectively. $\blacksquare$

\subsection{Harish-Chandra Schwartz spaces of Whittaker functions}\label{Section HCS spaces}

From this section and until the end of Section \ref{Section good sections Whittaker}, we assume that $G$ is quasi-split. Thus $P_0=B$ is a Borel subgroup, $T=M_0$ is a maximal torus and $A_0$ is the maximal split torus in $T$. Let $\xi: N_0(F)\to \mathbb{S}^1$ be a generic character. For every $\lambda\in \cA_0^*$, we denote by $\cC_{\lambda}(N_0(F)\backslash G(F),\xi)$ the space of functions $W:G(F)\to \C$ satisfying:
\begin{itemize}
\item $W(ug)=\xi(u)W(g)$ for every $(u,g)\in N_0(F)\times G(F)$;
\item If $F$ is $p$-adic: $W$ is right-invariant by a compact-open subgroup $J\subseteq G(F)$ and there exists $R>0$ such that for every $d>0$ we have an inequality
$$\displaystyle \lvert W(tk)\rvert \ll \left(\prod_{\alpha\in \Delta}(1+t^\alpha)^{R}\right)\Xi^G_\lambda(t)\overline{\sigma}(t)^{-d},\;\;\; t\in T(F),k\in K;$$
\item If $F$ is Archimedean: $W$ is smooth and there exists $R>0$ such that for every $u\in \cU(\g)$ and $d>0$ we have an inequality
$$\displaystyle \lvert (R(u)W)(tk)\rvert \ll \left(\prod_{\alpha\in \Delta}(1+t^\alpha)^{R}\right) \Xi^G_\lambda(t)\overline{\sigma}(t)^{-d},\;\;\; t\in T(F),k\in K.$$
\end{itemize}
Then, $\cC_{\lambda}(N_0(F)\backslash G(F),\xi)$ has a natural locally convex topology making it into a LF space. Notice that we have $\cC_{\lambda}(N_0(F)\backslash G(F),\xi)=\cC_{w\lambda}(N_0(F)\backslash G(F),\xi)$ for every $w\in W^G$ and $\lambda\in \cA_0^*$ (by Proposition \ref{prop 1 uniform bounds matrix coefficients}(i)). The next lemma shows that $\cC_{\lambda}(N_0(F)\backslash G(F),\xi)$ is actually a Fr\'echet space in the Archimedean case and a strict LF space in the $p$-adic case (in particular it is complete). Although it is probably well-known (see \cite[Proposition 3.1]{Jac} for a similar result) for the sake of completeness we provide a full proof.

\begin{lem}\label{lem 1 HCS spaces}
Let $\lambda\in \cA_0^*$. Then, for every $R>0$ and $d>0$ there exists a continuous semi-norm $p_{R,d}$ on $\cC_{\lambda}(N_0(F)\backslash G(F),\xi)$ such that
$$\displaystyle \lvert W(tk)\rvert \leqslant p_{R,d}(W)\left(\prod_{\alpha\in \Delta}(1+t^\alpha)^{-R}\right)\delta_0(t)^{1/2}t^{\lvert \lambda\rvert}\overline{\sigma}(t)^{-d}$$
for every $W\in \cC_{\lambda}(N_0(F)\backslash G(F),\xi)$, $t\in T(F)$ and $k\in K$.
\end{lem}

\noindent\ul{Proof}: By the uniform boundedness principle it suffices to show that for every $R>0$, $d>0$ and $W\in \cC_{\lambda}(N_0(F)\backslash G(F),\xi)$ we have
$$\displaystyle \lvert W(tk)\rvert \ll\left(\prod_{\alpha\in \Delta}(1+t^\alpha)^{-R}\right)\delta_0(t)^{1/2}t^{\lvert \lambda\rvert}\overline{\sigma}(t)^{-d},\;\;\; t\in T(F), k\in K.$$
Notice that by Proposition \ref{prop 1 uniform bounds matrix coefficients}(iv), there exist $R_0>0$ and $d_0>0$ such that
$$\displaystyle \Xi^G_\lambda(t)\ll \left(\prod_{\alpha\in \Delta}(1+t^\alpha)^{R_0}\right)\delta_0(t)^{1/2}t^{\lvert \lambda\rvert}\overline{\sigma}(t)^{d_0},\;\;\; t\in T(F).$$
Indeed, for every $t\in T(F)$ there exists $w\in W^G$ such that $wtw^{-1}\in T^+$ and since $w$ admits a lift in $K$, by Proposition \ref{prop 1 uniform bounds matrix coefficients}(iv) there exists $d_0>0$ such that
$$\displaystyle \Xi^G_\lambda(t)=\Xi^G_\lambda(wtw^{-1})\ll \delta_0(t)^{1/2}(wtw^{-1})^{\lvert \lambda\rvert}\overline{\sigma}(t)^{d_0}=t^{w^{-1}\lvert \lambda\rvert-\lvert \lambda\rvert}\delta_0(t)^{1/2}t^{\lvert \lambda\rvert}\overline{\sigma}(t)^{d_0}.$$
As $\lvert \lambda\rvert\in \overline{(\cA_0^*)^+}$, $w^{-1}\lvert \lambda\rvert-\lvert \lambda\rvert$ is a nonnegative linear combination of simple roots. Hence, there exists $R_0>0$ (which we can of course choose independently of $w$) so that
$$\displaystyle t^{w^{-1}\lvert \lambda\rvert-\lvert \lambda\rvert}\leqslant \prod_{\alpha\in \Delta}(1+t^\alpha)^{R_0}.$$
Therefore, we are left with showing that for every $R>0$, $d>0$ and $W\in \cC_{\lambda}(N_0(F)\backslash G(F),\xi)$ we have
$$\displaystyle \lvert W(tk)\rvert \ll\left(\prod_{\alpha\in \Delta}(1+t^\alpha)^{-R}\right)\Xi^G_\lambda(t)\overline{\sigma}(t)^{-d},\;\;\; t\in T(F), k\in K.$$
In the $p$-adic case we can prove the following stronger inequality by essentially the same argument as in the proof of \cite[Lemme 3.7]{Wa2}: for every compact-open subgroup $J$ of $G(F)$, there exists $c=c_J>0$ such that for every $W\in \cC_\lambda(N_0(F)\backslash G(F),\xi)^J$ and $d>0$ we have
$$\displaystyle \lvert W(tk)\rvert \ll \left(\prod_{\alpha\in \Delta}\mathbf{1}_{]0,c]}(t^\alpha)\right) \Xi^G_\lambda(t)\overline{\sigma}(t)^{-d},\;\;\; t\in T(F), k\in K$$
where $\mathbf{1}_{]0,c]}$ stands for the characteristic function of the interval $]0,c]$. Thus, we only consider the Archimedean case. Let $W\in \cC_{\lambda}(N_0(F)\backslash G(F),\xi)$ and $d>0$. Clearly, it suffices to show the existence of $R>0$ such that for every $\Delta$-tuple $\underline{N}=(N_\alpha)_{\alpha\in \Delta}$ of nonnegative integers we have
\begin{align}\label{eq 1 HCS Whittaker functions}
\displaystyle \left(\prod_{\alpha\in \Delta}t^{N_\alpha \alpha}\right)\lvert W(tk)\rvert \ll \left(\prod_{\alpha\in \Delta}(1+t^\alpha)^{R}\right) \Xi^G_\lambda(t)\overline{\sigma}(t)^{-d},\;\;\; t\in T(F), k\in K.
\end{align}
Let $d\xi:\n_0(F)\to \C$ denote the differential of $\xi$ at the origin. Since $\xi$ is generic for every $\alpha\in \Delta$ there exists $X_\alpha\in \n_0(F)$ such that $d\xi(X_\alpha)=1$ and $\Ad(a)X_\alpha=\alpha(a)X_\alpha$ for every $a\in A_0(F)$. We make such a choice. Fix a norm $\lVert .\rVert$ on $\n_0(F)$. As $T(F)/A_0(F)$ is compact, there exists $c>0$ such that $\lVert \Ad(tk)^{-1}X_\alpha\rVert \leqslant c t^{-\alpha}$ for every $t\in T(F)$, $k\in K$ and $\alpha\in \Delta$. Then, for every $\Delta$-tuple $\underline{N}=(N_\alpha)_{\alpha\in \Delta}$ of nonnegative integers, setting $u_{\underline{N}}=\prod_{\alpha\in \Delta}X_\alpha^{N_\alpha}\in \cU(\n_0)$ (the product being taken in some fixed order), we have
\[\begin{aligned}
\displaystyle \lvert W(tk)\rvert=\left\lvert \prod_{\alpha\in \Delta}d\xi(X_\alpha)^{N_\alpha} W(tk)\right\rvert & =\left\lvert(L(-\epsilon u_{\underline{N}})W)(tk)\right\rvert=\left\lvert\left(R(\epsilon\Ad(tk)^{-1}u_{\underline{N}}) W\right)(tk)\right\rvert \\
 & \leqslant \left(\prod_{\alpha\in \Delta}t^{-N_\alpha \alpha} \right) \sup_{u\in \cK_N} \lvert (R(u)W)(tk)\rvert
\end{aligned}\] 
for every $t\in T(F)$, $k\in K$ where $\epsilon$ is a certain sign and $\cK_N$ is the compact subset of $\cU(\n_0)$ consisting of products of $N=\sum_{\alpha\in \Delta}N_\alpha$ elements of $\n_0(F)$ of norm smaller than $c$. By definition of $\cC_{\lambda}(N_0(F)\backslash G(F),\xi)$, there exists $R>0$ such that for every integer $N\geqslant 0$ we have
$$\displaystyle \sup_{u\in \cK_N} \lvert (R(u)W)(tk)\rvert\ll \left(\prod_{\alpha\in \Delta}(1+t^\alpha)^{R}\right)\Xi^G_\lambda(t)\overline{\sigma}(t)^{-d},\;\;\; t\in T(F), k\in K.$$
This shows \ref{eq 1 HCS Whittaker functions} and ends the proof of the lemma. $\blacksquare$

\subsection{Uniform bounds for families of Whittaker functions}\label{Section bounds Whittaker}

We continue with the setting of the previous sections, still assuming that $G$ is quasi-split and fixing a generic character $\xi: N_0(F)\to \mathbb{S}^1$. Let $P=MN$ be a standard parabolic subgroup, $\sigma\in \Temp(M)$ and set $\pi_{\lambda}=i_P^G(\sigma_\lambda)$ for every $\lambda\in \cA_{M,\C}^*$. We identify as before the space of $\pi_\lambda$ with $\pi_K:=i_{K_P}^K(\sigma_{\mid K_M})$. We assume given a family $J_\lambda\in \Hom_{N_0}(\pi_\lambda,\xi)$ of Whittaker functionals on $\pi_\lambda=i_P^G(\sigma_{\lambda})$ for $\lambda\in \cA_{M,\C}^*$ i.e. a family of continuous linear forms $J_\lambda: \pi_\lambda\to \C$ satisfying $J_\lambda\circ \pi_\lambda(u)=\xi(u)J_\lambda$ for every $u\in N_0(F)$ and $\lambda\in \cA_{M,\C}^*$. We moreover suppose that the family $\lambda\mapsto J_\lambda\in (\pi_\lambda)'$ is holomorphic in the sense of Section \ref{Section representations}. By Frobenius reciprocity, for every $\lambda\in \cA_{M,\C}^*$, $J_\lambda$ induces a continuous $G(F)$-equivariant linear map
$$\displaystyle \widetilde{J_\lambda}: \pi_\lambda\to C^\infty(N_0(F)\backslash G(F),\xi)$$
where $C^\infty(N_0(F)\backslash G(F),\xi)$ stands for the space of all smooth functions $W:G(F)\to \C$ such that $W(ug)=\xi(u)W(g)$ for every $(u,g)\in N_0(F)\times G(F)$. Recall that in Section \ref{Section groups} we have defined a strict partial order $\prec$ on $\cA_0^*$.

\begin{prop}\label{prop 1 uniform bounds Whittaker functions}
\begin{enumerate}[(i)]
\item For every $\lambda\in \cA_{M,\C}^*$ and $\mu\in \cA_0^*$ such that $\lvert \Re(\lambda)\rvert\prec \mu$, the image of $\widetilde{J_\lambda}$ is included in $\cC_{\mu}(N_0(F)\backslash G(F),\xi)$ and the resulting linear map (that we will still denote by $\widetilde{J_\lambda}$)
$$\displaystyle \pi_\lambda\to \cC_{\mu}(N_0(F)\backslash G(F),\xi)$$
is continuous.

\item Let $\mu\in (\cA^G_0)^*$ and set $\cU[\prec \mu]=\{\lambda \in (\cA^G_{M,\C})^*\mid \lvert \Re(\lambda)\rvert\prec \mu\}$ (an open subset of $(\cA^G_{M,\C})^*$). Then, the family of continuous linear maps
$$\displaystyle \lambda\in \cU[\prec \mu]\mapsto \widetilde{J_\lambda}\in \Hom_{G(F)}(\pi_\lambda, \cC_{\mu}(N_0(F)\backslash G(F),\xi))$$
is analytic in the sense that for every analytic section $\lambda\mapsto e_\lambda\in \pi_\lambda$ the resulting map $\lambda\in \cU[\prec \mu]\mapsto \widetilde{J_\lambda}(e_\lambda)\in \cC_{\mu}(N_0(F)\backslash G(F),\xi)$ is analytic.
\end{enumerate}
\end{prop}

\noindent\ul{Proof}: We will show the following:
\begin{num}
\item\label{eq 1 uniform bounds Whittaker functions} For every compact subset $\cK\subset \cA_{M,\C}^*$, there exists $R>0$ and a continuous semi-norm $p$ on $\pi_K$ such that
$$\displaystyle \left\lvert J_\lambda(\pi_\lambda(t)e)\right\rvert\leqslant p(e)\Xi_{\Re(\lambda)}(t) \left(\prod_{\alpha\in \Delta}(1+t^\alpha)^{R}\right)$$
for every $\lambda\in \cK$, $t\in T(F)$ and $e\in \pi_K$.
\end{num}

Before proving \ref{eq 1 uniform bounds Whittaker functions}, we explain how the proposition can be deduced from this. The first part can readily be inferred from \ref{eq 1 uniform bounds Whittaker functions} together with the following inequality which is a consequence of Proposition \ref{prop 1 uniform bounds matrix coefficients} (iv) and the Cartan decomposition $G(F)=KT^+K$: If $\lvert \Re(\lambda)\rvert\prec \mu$ then for every $d>0$ we have
\begin{align}\label{eq 2 uniform bounds Whittaker functions}
\displaystyle \Xi_{\Re(\lambda)}(g)\ll \Xi_{\mu}(g)\overline{\sigma}(g)^{-d},\;\;\; g\in G(F).
\end{align}
For (ii), by the criterion \ref{eq 1 TVS} and since for every $e\in \pi_K$ and $g\in G(F)$ the map $\lambda\in \cA_{M,\C}^*\mapsto \widetilde{J_\lambda}(e)(g)=J_\lambda(\pi_\lambda(g)e)$ is analytic, it suffices to check that for every $e\in \pi_K$ the map $\lambda\in \cU[\prec \mu]\mapsto \widetilde{J_\lambda}(e)\in \cC_{\mu}(N_0(F)\backslash G(F),\xi)$ is continuous. Let $\lambda_0\in \cU[\prec \mu]$ and $e\in \pi_K$. Then, by definition of the topology on $\cC_{\mu}(N_0(F)\backslash G(F),\xi)$, we need to show the following: there exists $R>0$ such that for every $\epsilon>0$ and $d>0$ if $\lambda\in \cU[\prec \mu]$ is sufficiently close to $\lambda_0$ then
\begin{align}\label{eq 3 uniform bounds Whittaker functions}
\displaystyle \sup_{\substack{t\in T(F) \\ k\in K}} \Xi_\mu(t)^{-1}\left(\prod_{\alpha\in \Delta}(1+t^\alpha)^{-R}\right)\overline{\sigma}(t)^d \left\lvert J_\lambda(\pi_\lambda(tk)e)- J_{\lambda_0}(\pi_{\lambda_0}(tk)e)\right\rvert<\epsilon.
\end{align}
Let $R>0$ and $p$ be a continuous semi-norm on $\pi_K$ so that \ref{eq 1 uniform bounds Whittaker functions} is satisfied on a compact neighborhood $\cK\subseteq \cU[\prec \mu]$ of $\lambda_0$. Then, there exists $C>0$ such that
\[\begin{aligned}
\displaystyle \sup_{\substack{t\in T(F) \\ k\in K}} \Xi_\mu(t)^{-1}\left(\prod_{\alpha\in \Delta}(1+t^\alpha)^{-R}\right)\overline{\sigma}(t)^{d+1} \left\lvert J_\lambda(\pi_\lambda(tk)e)\right\rvert<C
\end{aligned}\]
for every $\lambda\in \cK$. Indeed, this follows from the (easily checked) fact that the inequality \ref{eq 2 uniform bounds Whittaker functions} can be made uniform on $\cK$. Set $M=2\epsilon^{-1}C$. Let $T[> M]$ denote the subset of elements $t\in T(F)$ such that $\overline{\sigma}(t)> M$ and $T[\leqslant M]$ be its complement. Then, the above inequality easily implies
\begin{align}\label{eq 4 uniform bounds Whittaker functions}
\displaystyle \sup_{\substack{t\in T[> M] \\ k\in K}} \Xi_\mu(t)^{-1}\left(\prod_{\alpha\in \Delta}(1+t^\alpha)^{-R}\right)\overline{\sigma}(t)^{d} \left\lvert J_\lambda(\pi_\lambda(tk)e)\right\rvert<\frac{\epsilon}{2}
\end{align}
for every $\lambda\in \cK$. On the other hand, $T[\leqslant M]$ is compact modulo $A_G(F)$ so that by continuity of the map $(\lambda,g)\mapsto J_\lambda(\pi_{\lambda}(g)e)$ for $\lambda\in (\cA_{M,\C}^G)^*$ sufficiently close to $\lambda_0$ we have
\[\begin{aligned}
\displaystyle \sup_{\substack{t\in T[\leqslant M] \\ k\in K}} \Xi_\mu(t)^{-1}\left(\prod_{\alpha\in \Delta}(1+t^\alpha)^{-R}\right)\overline{\sigma}(t)^d \left\lvert J_\lambda(\pi_\lambda(tk)e)- J_{\lambda_0}(\pi_{\lambda_0}(tk)e)\right\rvert<\epsilon.
\end{aligned}\]
Together with \ref{eq 4 uniform bounds Whittaker functions} this readily implies \ref{eq 3 uniform bounds Whittaker functions} thus showing (ii) provided \ref{eq 1 uniform bounds Whittaker functions} is satisfied.

We now prove \ref{eq 1 uniform bounds Whittaker functions}. Set $\widetilde{\pi}_\lambda=i_P^G(\widetilde{\sigma}_\lambda)$ that we identify with $\widetilde{\pi}_K:=i_{K_P}^K(\widetilde{\sigma}_{\mid K_M})$ for every $\lambda\in \cA_{M,\C}^*$. We fix a pairing $\langle.,.\rangle$ between $\pi_K$ and $\widetilde{\pi}_K$ as in Section \ref{Section matrix coefficients} that identifies $\widetilde{\pi}_{-\lambda}$ with the smooth contragredient of $\pi_{\lambda}$ for every $\lambda\in \cA_{M,\C}^*$. This pairing induces an embedding $\widetilde{\pi}_K\subseteq \pi_K'$.

First, we treat the $p$-adic case. Since $\xi$ is generic, we can choose for every $\alpha\in \Delta$ a vector $X_\alpha\in \n_0(F)$ so that $\xi(\exp(X_\alpha))\neq 1$ and $\Ad(a)X_\alpha=\alpha(a)X_\alpha$ for all $a\in A_0(F)$. Let $J\subseteq K$ be a compact-open subgroup. Then, there exists $C=C_J>0$ such that for every $\alpha\in \Delta$ and $t\in T(F)$ with $t^\alpha\geqslant C$ we have $t^{-1}\exp(X_\alpha)t\in J$. Let $e\in (\pi_K)^J$, $\lambda\in \cA_{M,\C}^*$ and $t\in T(F)$. Then, if there exists $\alpha\in \Delta$ such that $t^\alpha\geqslant C$ we have
\[\begin{aligned}
\displaystyle \xi(\exp(X_\alpha)) J_\lambda(\pi_{\lambda}(t)e)=J_\lambda(\pi_{\lambda}(t)\pi_{\lambda}(t^{-1}\exp(X_\alpha)t)e)=J_\lambda(\pi_{\lambda}(t)e)
\end{aligned}\]
hence (as $\xi(\exp(X_\alpha))\neq 1$) $J_\lambda(\pi_{\lambda}(t)e)=0$. Thus,
$$\displaystyle J_\lambda(\pi_{\lambda}(t)e)=0$$
unless $t^\alpha<C$ for every $\alpha\in \Delta$. On the other hand, there exists a compact-open subgroup $J'\subseteq K$ such that $J'\subseteq \Ker(\xi)tJt^{-1}$ as soon as $t^\alpha<C$ for every $\alpha\in \Delta$ where we have denoted by $\Ker(\xi)$ the kernel of $\xi$. Hence, if $t^\alpha<C$ for every $\alpha\in \Delta$ we have
$$\displaystyle J_\lambda(\pi_{\lambda}(t)e)=J_\lambda(\pi(e_{J'})\pi_{\lambda}(t)e)$$
where $e_{J'}:=\vol(J')^{-1}\mathbf{1}_{J'}$. The family of functionals $\lambda\mapsto J_\lambda\circ\pi(e_{J'})$ is represented by an analytic section $\lambda\in \cA_{M,\C}^*\mapsto e^\vee_\lambda\in \widetilde{\pi}_{\lambda}=i_P^G(\widetilde{\sigma}_\lambda)$ and by what we just saw for every $t\in T(F)$ we have
\[\begin{aligned}
\displaystyle J_\lambda(\pi_{\lambda}(t)e)=\left\{
    \begin{array}{ll}
        \langle \pi_{\lambda}(t)e, e^\vee_\lambda\rangle & \mbox{ if } t^\alpha<C \mbox{ for every } \alpha\in \Delta, \\
        0 & \mbox{ otherwise.}
    \end{array}
\right.
\end{aligned}\]
The inequality \ref{eq 1 uniform bounds Whittaker functions} (with $R=0$) is now a consequence of Proposition \ref{prop 2 uniform bounds matrix coefficients}.

Next, we treat the Archimedean case. For this, we need to introduce another set of notation. Fix a $M(F)$-continuous Hilbert norm $\widetilde{q}_0$ on $\widetilde{\sigma}$ which is $K_M$-invariant and let $\widetilde{p}_0$ be the $K$-invariant Hilbert norm on $\widetilde{\pi}_K$ defined by
$$\displaystyle \widetilde{p}_0(e^\vee)=\left(\int_K \widetilde{q}_0(e^\vee(k))^2 \right)^{1/2},\;\;\; e^\vee\in \widetilde{\pi}_K.$$
Then, for every $\lambda\in \cA_{M,\C}^*$, $\widetilde{p}_0$ induces a $G(F)$-continuous norm on $\widetilde{\pi}_{\lambda}$ through the identification $\widetilde{\pi}_{\lambda}\simeq \widetilde{\pi}_K$. Let $C_K$ be a Casimir operator for $K$. Set, as in Lemma \ref{lem 1 representations}, $\widetilde{p}_\ell=\widetilde{p}_0\circ \widetilde{\pi}_K(1+C_K)^\ell$ for every integer $\ell\in \bZ$ and $\widetilde{\pi}_K^{(\ell)}$ for the completion of $\widetilde{\pi}_K$ with respect to $\widetilde{p}_\ell$ (a Banach space). Then, by Lemma \ref{lem 1 representations}, we have natural inclusions $\widetilde{\pi}_K^{(\ell)}\subseteq \pi_K'$ for every $\ell$, every equicontinuous subset of $\pi_K'$ is contained and bounded in some $\widetilde{\pi}_K^{(\ell)}$ where $\ell\leqslant 0$ and the family of norms $(\widetilde{p}_\ell)_{\ell\geqslant 0}$ generates the topology on $\widetilde{\pi}_K$. Let $\ell\leqslant 0$ and $m\geqslant -\ell$ be integers and $\varphi\in C_c^{2m}(G(F))$ then we have
\begin{align}\label{eq 4bis uniform bounds Whittaker functions}
\displaystyle \pi'_{\lambda}(\varphi)\Lambda\in \widetilde{\pi}_K^{(\ell+m)}
\end{align}
for every $\lambda\in \cA_{M,\C}^*$ and $\Lambda\in \widetilde{\pi}_K^{(\ell)}$. Indeed, we have
\[\begin{aligned}
\displaystyle \widetilde{p}_{\ell+m}(\widetilde{\pi}_\lambda(\varphi)e^\vee) & =\widetilde{p}_0\left(\widetilde{\pi}_\lambda(L(1+C_K)^{\ell+m}R(1+C_K)^{-\ell}\varphi)\widetilde{\pi}_K(1+C_K)^\ell e^\vee \right) \\
 & \leqslant \left\lVert L(1+C_K)^{\ell+m}R(1+C_K)^{-\ell}\varphi\right\rVert_{L^1} \times \sup_{g\in \Supp(\varphi)} \lVert \widetilde{\pi}_\lambda(g)\rVert_{\widetilde{p}_0} \times \widetilde{p}_\ell(e^\vee)  
\end{aligned}\]
for every $\lambda\in \cA_{M,\C}^*$ and $e^\vee\in \widetilde{\pi}_K$ where $\left\lVert L(1+C_K)^{\ell+m}R(1+C_K)^{-\ell}\varphi\right\rVert_{L^1}$ denotes the $L^1$-norm of $L(1+C_K)^{\ell+m}R(1+C_K)^{-\ell}\varphi\in C_c(G(F))$ and for every $g\in G(F)$ we have denoted by $\lVert \widetilde{\pi}_\lambda(g)\rVert_{\widetilde{p}_0}$ the operator norm of $\widetilde{\pi}_{\lambda}(g)$ with respect to the norm $\widetilde{p}_0$. This proves \ref{eq 4bis uniform bounds Whittaker functions} by density of $\widetilde{\pi}_K$ in $\widetilde{\pi}_K^{(\ell')}$ for every $\ell'$. Moreover, since the function $\lambda\in \cA_{M,\C}^*\mapsto \sup_{g\in \Supp(\varphi)} \lVert \widetilde{\pi}_\lambda(g)\rVert_{\widetilde{p}_0}$ is easily seen to be locally bounded, we have that the operator norm of $\pi'_\lambda(\varphi)$ seen as a continuous linear map from $\widetilde{\pi}_K^{(\ell)}$ to $\widetilde{\pi}_K^{(\ell+m)}$ is locally bounded in $\lambda$.

Let $\cK\subseteq \cA_{M,\C}^*$ be a compact subset. By continuity of $\lambda\mapsto J_\lambda\in \pi_K'$ for the weak-$\star$ topology and the uniform boundedness principle, the family $\{J_\lambda\mid \lambda\in \cK \}$ is equicontinuous. Hence, there exists $\ell\leqslant 0$ such that this family is included and bounded in $\widetilde{\pi}_K^{(\ell)}$ and, by what we have just seen, for every $m\geqslant -\ell$ and $\varphi\in C_c^{2m}(G(F))$ the family
$$\displaystyle \left\{J_\lambda\circ \pi_{\lambda}(\varphi) \mid \lambda\in \cK \right\}$$
is included and bounded in $\widetilde{\pi}_K^{(\ell+m)}$. Let $p$ and $\widetilde{p}$ be continuous semi-norms on $\pi_K$ and $\widetilde{\pi}_K$ as in Proposition \ref{prop 2 uniform bounds matrix coefficients}. Then, for $m$ sufficiently large $\widetilde{p}$ extends (uniquely) to a continuous semi-norm on $\widetilde{\pi}_K^{(m+\ell)}$ and the inequality of Proposition \ref{prop 2 uniform bounds matrix coefficients} still holds for every $(e,e^\vee)\in \pi_K\times \widetilde{\pi}_K^{(m+\ell)}$ (by density of  $\widetilde{\pi}_K$ in $\widetilde{\pi}_K^{(m+\ell)}$). Therefore, we have obtained:
\begin{num}
\item\label{eq 5 uniform bounds Whittaker functions} For $m$ sufficiently large and every $\varphi\in C_c^{2m}(G(F))$ there exists a constant $C>0$ such that
$$\displaystyle \lvert J_\lambda(\pi_{\lambda}(\varphi)\pi_{\lambda}(g)e)\rvert \leqslant C p(e) \Xi^G_{\Re(\lambda)}(g)$$
for all $\lambda\in \cK$, $g\in G(F)$ and $e\in \pi_K$.
\end{num}

Let $\overline{B}=T\overline{N}_0$ be the Borel subgroup opposite to $B$ and let $Y_1,\ldots,Y_b$ be an $\bR$-basis of $\overline{\fb}(F)$. Set $\Delta_{\overline{B}}=Y_1^2+\ldots+Y_b^2\in \cU(\overline{\fb})$ and let $m$ be a positive integer that we assume sufficiently large in what follows. By elliptic regularity (\cite[Lemma 3.7]{BK}) there exist $\varphi_{\overline{B}}^1\in C_c^{m'}(\overline{B}(F))$ and $\varphi_{\overline{B}}^2\in C_c^\infty(\overline{B}(F))$ where $m'=2m-\dim(\overline{B})-1$ such that
$$\displaystyle \varphi_{\overline{B}}^1 \ast \Delta_{\overline{B}}^m+\varphi_{\overline{B}}^2=\delta_1^{\overline{B}}$$
in the sense of distributions where $\delta_1^{\overline{B}}$ denotes the Dirac distribution at $1$ on $\overline{B}(F)$. Let $\varphi_N\in C_c^\infty(N_0(F))$ be such that
$$\displaystyle \int_{N_0(F)} \varphi_N(u) \xi(u) du=1.$$
Then, setting $\varphi^i=\varphi_N\ast \varphi_{\overline{B}}^i$ for $i=1,2$, for every $\lambda\in \cA_{M,\C}^*$, $e\in \pi_K$ and $t\in T(F)$, we have
\begin{align}\label{eq 6 uniform bounds Whittaker functions}
\displaystyle J_\lambda(\pi_{\lambda}(t)e) & =J_\lambda(\pi_\lambda(\varphi_{\overline{B}}^1)\pi_\lambda(\Delta_{\overline{B}}^m)\pi_\lambda(t)e)+J_\lambda(\pi_\lambda(\varphi_{\overline{B}}^2)\pi_\lambda(t)e) \\
\nonumber & =J_\lambda(\pi_\lambda(\varphi^1)\pi_\lambda(\Delta_{\overline{B}}^m)\pi_\lambda(t)e)+J_\lambda(\pi_\lambda(\varphi^2)\pi_\lambda(t)e).
\end{align}
Noticing that $\varphi^i\in C_c^{m'}(G(F))$ for $i=1,2$, from \ref{eq 5 uniform bounds Whittaker functions} we deduce that for $m$ sufficiently large there exists $C>0$ such that
\begin{align}\label{eq 7 uniform bounds Whittaker functions}
\displaystyle \lvert J_\lambda(\pi_{\lambda}(t)e)\rvert \leqslant C(p(\pi_{\lambda}(\Ad(t)^{-1}\Delta_{\overline{B}}^m)e)+p(e))\Xi^G_{\Re(\lambda)}(t)
\end{align}
for every $\lambda\in \cK$, $t\in T(F)$ and $e\in \pi_K$. To get \ref{eq 1 uniform bounds Whittaker functions} it only remains to notice that there exist $R>0$ and a continuous semi-norm $q$ on $\pi_K$ such that
$$\displaystyle p(\pi_\lambda(\Ad(t)^{-1}\Delta_{\overline{B}}^m)e)\leqslant \left(\prod_{\alpha\in \Delta}(1+t^\alpha)^R \right)q(e)$$
for every $\lambda\in \cK$, $t\in T(F)$ and $e\in \pi_K$. This follows from the fact that the function $(\lambda,t)\in \cA_{M,\C}^*\times T(F)\mapsto \pi_\lambda(\Ad(t)^{-1}\Delta_{\overline{B}}^m)$ is polynomial (see \cite[Proposition 1]{Del}). $\blacksquare$

\subsection{Application to the existence of good sections for Whittaker models}\label{Section good sections Whittaker}

We continue with the setting of the previous section: $G$ is quasi-split, $\xi:N_0(F)\to \mathbb{S}^1$ is a generic character, $P=MN$ is a standard parabolic subgroup, $\sigma\in \Temp(M)$ and we set $\pi_{\lambda}=i_P^G(\sigma_\lambda)$ for every $\lambda\in \cA_{M,\C}^*$ whose space is identified as before with $\pi_K:=i_{K_P}^K(\sigma_{\mid K_M})$. Let $w_0\in W^G$ be the longest element of the Weyl group and choose a lift $\widetilde{w_0}\in G(F)$. Let $\xi^-$ be the generic character of $\overline{N}_0(F)=w_0N_0(F)w_0^{-1}$ defined by $\xi^-(\overline{u})=\xi(\widetilde{w_0}\overline{u}\widetilde{w_0}^{-1})$, $\overline{u}\in \overline{N}_0(F)$. Assume that $\sigma$ is generic with respect to the restriction of $\xi^-$ to $\overline{N}_0(F)\cap M(F)$ i.e. there exists a continuous nonzero linear form $\ell:\sigma\to \C$ such that $\ell\circ \sigma(\overline{u})=\xi^-(\overline{u})\ell$ for every $\overline{u}\in \overline{N}_0(F)\cap M(F)$. Then, the construction of Jacquet's integral (\cite[Sect. 15.4]{Wall2}, \cite{CS}) provides us with a holomorphic family of Whittaker functionals
$$\displaystyle \lambda\in \cA_{M,\C}^*\mapsto J_\lambda\in \Hom_{N_0}(\pi_\lambda,\xi)$$
as in the previous section which is everywhere non-vanishing. For every $\lambda\in \cA_{M,\C}^*$ we denote by $\cW(\pi_{\lambda},\xi)$ the corresponding Whittaker model i.e. the space of functions of the form $g\in G(F)\mapsto J_\lambda(\pi_{\lambda}(g)e)$ for $e\in \pi_K$.

\begin{cor}\label{cor good sections Whittaker models}
For every $\lambda_0\in (\cA^G_{M,\C})^*$ and $W_0\in\cW(\pi_{\lambda_0},\xi)$ there exists a map
$$\displaystyle \lambda\in (\cA^G_{M,\C})^*\mapsto W_\lambda\in \cW(\pi_{\lambda},\xi)$$
such that:
\begin{itemize}
\item for every $\mu\in (\cA^G_0)^*$ and $\lambda\in\cU[\prec \mu]=\{\lambda \in (\cA^G_{M,\C})^*\mid \lvert \Re(\lambda)\rvert\prec \mu\}$ we have $W_\lambda\in \cC_\mu(N_0(F)\backslash G(F),\xi)$ and the resulting map
$$\displaystyle \lambda\in \cU[\prec \mu]\mapsto W_\lambda\in \cC_\mu(N_0(F)\backslash G(F),\xi)$$
is analytic;
\item $W_{\lambda_0}=W_0$.
\end{itemize}
\end{cor}

\noindent\ul{Proof}: Indeed, there exists $e\in \pi_K$ such that (with the notation of Section \ref{Section bounds Whittaker}), $W_0=\widetilde{J}_{\lambda_0}(e)$ and then it suffices to set $W_\lambda=\widetilde{J}_{\lambda}(e)$ for every $\lambda\in (\cA^G_{M,\C})^*$: the required properties immediately follows from Proposition \ref{prop 1 uniform bounds Whittaker functions}. $\blacksquare$

\subsection{Holomorphic continuation of certain functions}\label{Section automatic holomorphic cont}

For every $C,D\in \bR\cup\{-\infty \}$ with $D>C$ we set $\cH_{>C}=\{s\in \C\mid \Re(s)>C \}$ and $\cH_{]C,D[}=\{s\in \C\mid C<\Re(s)<D  \}$. A {\em vertical strip} is a subset of $\C$ which is the closure of $\cH_{]C,D[}$ for some $C,D\in \bR$ with $D>C$. Let $M$ be a complex analytic manifold and $C\in \bR\cup\{-\infty \}$. Then, we say that a holomorphic function $Z:\cH_{>C}\times M\to \C$ is {\em of finite order in vertical strips in the first variable locally uniformly in the second variable} if there exists $d\in \bR$ such that for every vertical strip $V\subseteq \cH_{>C}$ and every compact subset $\cK_M\subseteq M$, we have
$$\displaystyle \sup_{(s,t)\in V\times \cK_M}e^{-\lvert s\rvert^{d}}\left\lvert Z(s,t)\right\rvert<\infty.$$
If the above inequality holds with $d=0$ for every vertical strip $V\subseteq \cH_{>C}$ and every compact subset $\cK_M\subseteq M$, we also say that $Z$ is {\em bounded in vertical strips in the first variable locally uniformly in the second variable}. Similarly, we say that $Z$ is {\em rapidly decreasing in vertical strips in the first variable locally uniformly in the second variable} if for every vertical strip $V\subseteq \cH_{>C}$, every $N>0$ and every compact subset $\cK_M\subseteq M$ we have
$$\displaystyle \sup_{(s,t)\in V\times \cK_M}(1+\lvert s\rvert)^N\left\lvert Z(s,t)\right\rvert<\infty.$$
When $M$ is a point, i.e. $Z$ is just a function of one complex variable, we will just say ``of finite order in vertical strips'', resp. ``bounded in vertical strips'', resp. ``rapidly decreasing in vertical strips''. The following proposition will be needed to extend Theorem \ref{theo 1 intro} from nearly tempered representations to all (generic) representations.

\begin{prop}\label{prop 1 holomorphic continuation}
Let $M$ be a connected complex analytic variety, $U\subseteq M$ a nonempty connected relatively compact open subset and
$$\displaystyle Z_+,Z_-: \C\times U\to \C$$
be two holomorphic functions satisfying
\begin{align}\label{eq 1 holo contn}
\displaystyle Z_+(s,t)=Z_-(-s,t)
\end{align}
for every $(s,t)\in \C\times U$. Moreover, we assume that for each $\epsilon\in \{\pm\}$ the following two conditions are satisfied
\begin{enumerate}[(1)]
\item $Z_\epsilon$ is of finite order in vertical strips in the first variable locally uniformly in the second variable;
\item For every connected relatively compact open subset $U'\subseteq M$ containing $U$, there exists $C>0$ such that $Z_\epsilon$ admits a (necessarily unique) holomorphic continuation to $\cH_{>C}\times U'$ which is again of finite order in vertical strips in the first variable locally uniformly in the second variable.
\end{enumerate}
Then, $Z_+$  and $Z_-$ extend to holomorphic functions on all of $\C\times M$ which are of finite order in vertical strips in the first variable locally uniformly in the second variable and still satisfy \ref{eq 1 holo contn}.
\end{prop}

\noindent\ul{Proof}: Let $D>0$ and $U'\subseteq M$ be a connected relatively compact open subset containing $U$. Clearly, we just need to show that $Z_+$, $Z_-$ extend to holomorphic functions on $\cH_{]-D,D[}\times U'$ of finite order in vertical strips in the first variable locally uniformly in the second variable (the functional equation \ref{eq 1 holo contn} will then hold automatically for these extensions by connectedness). Choosing another connected relatively compact open subset $U''\subseteq M$ containing the closure of $U'$ and using condition (2) for $U''$, we see that there exists $C>0$ and a nonnegative integer $n$ such that $(s,t)\mapsto e^{s^{n}}Z_{\pm}(s,t)$ extend to holomorphic functions on $\cH_{>C}\times U'\cup \C\times U$ which are rapidly decreasing in vertical strips in the first variable locally uniformly in the second variable. Moreover, up to multiplying $Z_+$ and $Z_-$ by $(s,t)\mapsto e^{s^{n}}$, we may assume that $n=0$. Then, for $D>C$ and $\epsilon\in \{\pm \}$, we set
$$\displaystyle \cZ_\epsilon^{<D}(s,t)=\frac{1}{2\pi}\left(\int_{-\infty}^{+\infty} \frac{Z_\epsilon(D+ix,t)}{D+ix-s}dx-\int_{-\infty}^{+\infty} \frac{Z_{-\epsilon}(D+ix,t)}{D+ix+s}dx \right)$$
for every $(s,t)\in \cH_{]-D,D[}\times U'$. The above expression converges absolutely and defines a holomorphic function on $\cH_{]-D,D[}\times U'$ which is bounded in vertical strips in the first variable locally uniformly in the first variable. Moreover, by the functional equation \ref{eq 1 holo contn} and since $Z_+$, $Z_-$ are holomorphic functions on $\C\times U$ that are rapidly decreasing in vertical strips in the first variable, by Cauchy's integration formula $\cZ_\epsilon^{<D}$ and $Z_\epsilon$ coincide on $\cH_{]-D,D[}\times U$ for $\epsilon\in\{\pm\}$. Therefore, $\cZ_\epsilon^{<D}$ is the sought-after extension. $\blacksquare$

\section{Zeta integrals and statement of the main theorems}\label{Part II}

\subsection{Notation}\label{section notation}

Let $F$ be a local field of characteristic zero and $E$ be either a quadratic extension of $F$ (inert case) or $F\times F$ (split case). We write $\lvert. \rvert_F$ and $\lvert .\rvert_E$ for the normalized absolute values of $F$ and $E$ respectively. Thus, in the split case we have $\lvert(\lambda,\mu)\rvert_E=\lvert \lambda\rvert_F \lvert \mu\rvert_F$ for every $(\lambda,\mu)\in E$ and in both cases we have $\lvert x\rvert_E=\lvert x\rvert_F^2$ for every $x\in F$. In the non-Archimedean case, we let $\cO_F$ and $\cO_E$ be the rings of integers of $F$ and $E$ respectively. We fix non-trivial additive characters $\psi': F\to \mathbb{S}^1$ and $\psi: E\to \mathbb{S}^1$ and we assume that {\em $\psi$ is trivial on $F$}. We denote by $\tau$ the unique element of $E$ such that $\psi(z)=\psi'(\Tra_{E/F}(\tau z))$ for every $z\in E$ where $\Tra_{E/F}$ stands for the trace of the extension $E/F$.

Let $n\geqslant 1$ be an integer. We write $G_n$ for the algebraic group $\GL_n$, $\mathfrak{gl}_n$ for its Lie algebra and we denote by $Z_n$, $B_n$, $A_n$ and $N_n$ the subgroups of scalar, resp. upper triangular, resp. diagonal, resp. unipotent upper triangular matrices in $G_n$. By a {\em standard} parabolic subgroup $P$ of $G_n$ we mean a parabolic subgroup containing $B_n$. We then write $P=MU$ for its unique Levi decomposition with $A_n\subseteq M$. We denote by $\delta_n$ and $\delta_{n,E}$ the modular characters of $B_n(F)$ and $B_n(E)$ respectively. We will denote by $K_n$ the standard maximal compact subgroup of $G_n(F)$ that is $K_n=G_n(\cO_F)$ in the $p$-adic case and $K_n=O(n)$ or $U(n)$ in the cases where $F=\bR$ or $\C$ respectively. For $g\in G_n$, we write ${}^t g$ for its transpose and $g_{i,j}$, $1\leqslant i,j\leqslant n$, for its entries. Also, for $a\in A_n$ we simply write $a_i$ for $a_{i,i}$ ($1\leqslant i\leqslant n$). Let $\cA^*=\cA_{A_n}^*$, $\cA_\C^*=\cA_{A_n,\C}^*$ and $(\cA^{G_n})^*=(\cA^{G_n}_{A_n})^*$. We identify $\cA^*$ with $\bR^n$ through the choice of the basis $\chi_1,\ldots,\chi_n$ of $X^*(A_n)$ where $\chi_i$ ($1\leqslant i\leqslant n$) is defined by $\chi_i(a)=a_i$ for every $a\in A_n$. Then, the closed negative Weyl chamber $\overline{(\cA^*)^+}$ with respect to $B_n$ is the cone of $n$-tuples $\lambda=(\lambda_1,\ldots,\lambda_n)\in \bR^n$ such that $\lambda_1\leqslant \ldots\leqslant \lambda_n$ and the subspace $(\cA^{G_n})^*$ consists of vectors $\lambda\in \bR^n$ with $\lambda_1+\ldots+\lambda_n=0$. Let $\lambda=(\lambda_1,\ldots,\lambda_n)\in \cA^*$. Then, defining $\lvert \lambda\rvert$ as in Section \ref{Section groups}, we have $\lvert \lambda\rvert=(\lambda_{w(1)},\ldots,\lambda_{w(n)})$ where $w$ is any permutation of $\{1,\ldots,n \}$ such that $\lambda_{w(1)}\leqslant\ldots\leqslant \lambda_{w(n)}$. Moreover, the character $a\mapsto a^\lambda$ of $A_n(E)$ is given by
$$\displaystyle a^\lambda=\lvert a_1\rvert_E^{\lambda_1}\ldots \lvert a_n\rvert_E^{\lambda_n}.$$
We set $\min(\lambda):=\min(\lambda_1,\ldots,\lambda_n)$. Notice that if $\lambda \in \overline{(\cA^*)^+}$ then $\min(\lambda)=\lambda_1$.

For $\pi\in \Irr(G_n(E))$ we will denote by $\omega_\pi$ its central character seen as a character on $E^\times$ through the natural identification $Z_n(E)\simeq E^\times$. Let $P=MU$ be a standard parabolic subgroup of $G_n$. Then, $M$ is of the form
$$\displaystyle M=G_{n_1}\times\ldots \times G_{n_k}$$
for some integers $n_1,\ldots,n_k$ with $n_1+\ldots+n_k=n$. Let $\tau_i\in \Irr(G_{n_i}(E))$, for every $1\leqslant i\leqslant k$, be an irreducible representation of $G_{n_i}(E)$ so that $\sigma=\tau_1\boxtimes \ldots\boxtimes \tau_k$ is an irreducible representation of $M(E)$. We will denote by
$$\displaystyle \tau_1\times\ldots\times \tau_k$$
the (normalized) induced representation $i_{P(E)}^{G_n(E)}(\sigma)$. A representation $\pi\in \Irr(G_n(E))$ is {\em generic} if it admits a nonzero Whittaker functional with respect to any (or equivalently one) generic character of $N_n(E)$. We will denote by $\Irr_{\gen}(G_n(E))$ the subset of generic representations in $\Irr(G_n(E))$. By \cite[Theorem 9.7]{Zel} and \cite[Theorem 6.2f]{Vog}, every $\pi\in \Irr_{\gen}(G_n(E))$ is isomorphic to a representation of the form $\tau_1\times\ldots \tau_k$ where for each $1\leqslant i\leqslant k$, $\tau_i$ is an essentially square-integrable (i.e. an unramified twist of a square-integrable) representation of some $G_{n_i}(E)$. We will say that a representation $\pi$ of $G_n(E)$ is {\em nearly tempered} if it is isomorphic to an induced representation of the form
$$\displaystyle (\tau_1\otimes \lvert \det \rvert_E^{\lambda_1})\times \ldots\times (\tau_k\otimes \lvert \det \rvert_E^{\lambda_k})$$
where for each $1\leqslant i\leqslant k$, $\tau_i\in \Pi_2(G_{n_i}(E))$ for some $n_i\geqslant 1$ and $\lambda_i$ is a real number with $\lvert \lambda_i\rvert<\frac{1}{4}$. Again by \cite[Theorem 9.7]{Zel} and also \cite{VS}, every nearly tempered representation is irreducible and generic. We will denote by $\Irr_{\ntemp}(G_n(E))\subseteq \Irr_{\gen}(G_n(E))$ the subset of nearly tempered representations.

We equip $\Temp(G_n(E))$ with a topology as follows. For each standard parabolic subgroup $P=MU$ of $G_n$ and $\sigma\in \Pi_2(M(E))$ the map $\lambda\in i\cA_M^*\mapsto i_{P(E)}^{G_n(E)}(\sigma_\lambda) \in \Temp(G_n(E))$ identifies a certain quotient of $i\cA_M^*$ with a subset of $\Temp(G_n(E))$. Then, we endow $\Temp(G_n(E))$ with the unique topology whose connected components are precisely these subsets equipped with the quotient topology from $i\cA_M^*$.

We define generic characters $\psi_n: N_n(E)\to \mathbb{S}^1$, $\psi'_n:N_n(F)\to \mathbb{S}^1$ by
$$\displaystyle \psi'_n(u)=\psi'\left((-1)^n\sum_{i=1}^{n-1} u_{i,i+1}\right) \mbox{ and } \psi_n(u)=\psi\left((-1)^n\sum_{i=1}^{n-1} u_{i,i+1}\right).$$
If $\pi$ is an irreducible generic representation of $G_n(E)$, we write $\cW(\pi,\psi_n)$ for its Whittaker model (with respect to $\psi_n$). For every Whittaker function $W\in C^\infty(N_n(E)\backslash G_n(E),\psi_n)$, we define $\widetilde{W}\in C^\infty(N_n(E)\backslash G_n(E),\psi_n^{-1})$ by
$$\displaystyle \widetilde{W}(g)=W(w_n{}^tg^{-1})\;\;\; (g\in G_n(E))$$
where $w_n=\begin{pmatrix} & & 1 \\ & \iddots & \\ 1 & & \end{pmatrix}$. Then, if $\pi$ is an irreducible generic representation of $G_n(E)$ the map $W\mapsto \widetilde{W}$ induces a (topological) isomorphism $\cW(\pi,\psi_n)\simeq \cW(\widetilde{\pi},\psi_n^{-1})$.

We let $\cS(F^n)=C_c^\infty(F^n)$ be the space of all locally constant and compactly supported complex-valued functions on $F^n$ in the $p$-adic case and $\cS(F^n)$ be the usual Schwartz space on $F^n$ in the Archimedean case. We denote by $\phi\mapsto \widehat{\phi}$ be the usual Fourier transform on $F^n$ defined using the additive character $\psi'$ and the corresponding autodual measure i.e. for every $\phi\in \cS(F^n)$ we have
$$\displaystyle \widehat{\phi}(x_1,\ldots,x_n)=\int_{F^n} \phi(y_1,\ldots,y_n)\psi'(x_1y_1+\ldots+x_ny_n)dy_1\ldots dy_n,\;\; (x_1,\ldots,x_n)\in F^n,$$
where the measure of integration is chosen such that $\widehat{\widehat{\phi}}(v)=\phi(-v)$. Finally, we write $e_n$ for the vector $e_n=(0,\ldots,0,1)\in F^n$.

\subsection{Asai $L$-functions and epsilon factors}\label{Section Asai factors}

Let $W_F$ be the Weil group of $F$ and set
$$\displaystyle W'_F=\left\{
    \begin{array}{ll}
        W_F\times SL_2(\C) & \mbox{ if } F \mbox{ is } p-\mbox{adic}, \\
        W_F & \mbox{ if } F \mbox{ is  Archimedean},
    \end{array}
\right.
$$
for the Weil-Deligne group of $F$. In the inert case, we define similarly the Weil-Deligne group $W'_E$ of $E$. An {\em admissible} complex representation of $W'_F$ (resp. $W'_E$) is by definition a continuous morphism $\phi:W'_F\to \GL(M)$ (resp. $\phi:W'_E\to GL(M)$), where $M$ is a finite dimensional complex vector space, which is semi-simple and algebraic when restricted to $SL_2(\C)$ (in the $p$-adic case). To any admissible complex representation $\phi:W'_F\to \GL(M)$, we associate a local $L$-factor $L(s,\phi)$ and a local $\epsilon$-factor $\epsilon(s,\phi,\psi')$ as in \cite[\S 3]{Ta} and \cite[\S 2.2]{GR}. In the $p$-adic case, $L(s,\phi)$ is of the form $P(q^{-s})$ where $P\in \C[T]$ is such that $P(0)=1$ whereas $\epsilon(s,\phi,\psi')$ is of the form $cq^{n(s-1/2)}$ where $n\in \bZ$ and $c=\epsilon(1/2,\phi,\psi')\in \C^\times$. In the Archimedean case, $L(s,\phi)$ is a product of $\dim(M)$ factors of the form $\pi^{-(s+s_0)/2}\Gamma((s+s_0)/2)$ for some $s_0\in \C$ where $\Gamma$ denotes the usual gamma function and $\epsilon(s,\varphi,\psi')$ is of the form $c Q^{s-1/2}$ where $Q\in \R_+^*$ and $c=\epsilon(1/2,\varphi,\psi')\in \C^\times$. Using these invariants, we define the local $\gamma$-factor of $\phi$ as
$$\displaystyle \gamma(s,\phi,\psi')=\epsilon(s,\phi,\psi')\frac{L(1-s,\widetilde{\phi})}{L(s,\phi)}$$
where $\widetilde{\phi}$ stands for the contragredient of $\phi$. Setting $\psi'_\lambda(x)=\psi'(\lambda x)$ for every $(\lambda,x)\in F^\times \times F$, we have
$$\displaystyle \epsilon(s,\phi,\psi'_\lambda)=(\det \phi)(\lambda)\lvert \lambda\rvert_F^{\dim(M)(s-\frac{1}{2})}\epsilon(s,\phi,\psi')$$
for every $\lambda\in F^\times$ where we have identified $\det \phi$ with a character of $F^\times$ through class field theory. Moreover, if we identify the absolute value $\lvert .\rvert_F$ to a character $W'_F\to \bR_+^*$ (again through class field theory), we have $L(s,\phi\otimes \lvert .\rvert_F^{s_0})=L(s+s_0,\phi)$ and $\epsilon(s,\phi\otimes \lvert .\rvert_F^{s_0},\psi')=\epsilon(s+s_0,\phi,\psi')$ for every $s_0\in \C$. Let $\eta_{E/F}$ be the quadratic character of $W'_F$ associated to the extension $E/F$ and set
$$\displaystyle \lambda_{E/F}(\psi')=\epsilon(\frac{1}{2},\eta_{E/F},\psi').$$
This is sometimes called the {\em Langlands constant} of the extension $E/F$. It is a fourth root of unity which is trivial in the split case. Moreover, in the inert case if $\phi$ is an admissible complex representation of $W'_E$ of dimension $n$, by the inductivity of $\epsilon$-factors in degree $0$ (\cite[Theorem 3.4.1]{Ta}), we have
\begin{align}\label{eq 1 As L functions}
\displaystyle \epsilon(s,\phi,\psi'_E)=\lambda_{E/F}(\psi')^{-n}\epsilon(s,\Ind_{W'_E}^{W'_F}(\phi),\psi')
\end{align}
where we have set $\psi'_E:=\psi'\circ\Tra_{E/F}$ and $\Ind_{W'_E}^{W'_F}$ denotes the functor of induction from $W'_E$ to $W'_F$ (which sends admissible representations to admissible representations).

\vspace{2mm}

Assume that we are in the inert case. Fix $s\in W'_F\setminus W'_E$. Let $\phi: W'_E\to \GL(M)$ be an admissible representation. We defined $\As(\phi):W'_F\to \GL(M\otimes M)$ by $\As(\phi)(w)=\phi(w)\otimes \phi(sws^{-1})$ for $w\in W'_E$ and $\As(\phi)(s)=(\Id_M\otimes \phi(s^2))\circ \iota$ where $\iota$ is the linear automorphism of $M\otimes M$ sending $u\otimes v$ to $v\otimes u$. Then, $\As(\phi)$ is an admissible representation of $W'_F$ and we set $L(s,\phi,\As)=L(s,\As(\phi))$, $\epsilon(s,\phi,\As,\psi')=\epsilon(s,\As(\phi),\psi')$, $\gamma(s,\phi,\As,\psi')=\gamma(s,\As(\phi),\psi')$ and call them the Asai $L$-function, $\epsilon$-factor and $\gamma$-factor of $\phi$ respectively. If $\phi$ decomposes as a direct sum
$$\displaystyle \phi=\phi_1\oplus \ldots\oplus \phi_k$$
of admissible representations of $W'_E$ then we have
\begin{align}\label{eq 2 As L functions}
\displaystyle \As(\phi)=\bigoplus_{i=1}^k \As(\phi_i)\oplus \bigoplus_{1\leqslant i<j\leqslant k}\Ind_{W'_E}^{W'_F}(\phi_i\otimes \phi_j^c)
\end{align}
where $\phi_j^c$ stands for the representation $w\mapsto \phi_j(sws^{-1})$. We have $\As(\phi\otimes \lvert .\rvert_E^{s_0})=\As(\phi)\otimes \lvert .\rvert_F^{2s_0}$ for every $s_0\in \C$. Moreover, if $\phi=\chi$ is a character of $W'_E$, that we identify to a character of $E^\times$ through class field theory, then $\As(\chi)$ is the character of $W'_F$ corresponding to the restriction of $\chi$ to $F^\times$.

Let $\pi\in \Irr(G_n(E))$. Then, the local Langlands correspondence for $\GL_n$ (\cite{La}, \cite{HT}, \cite{Hen}, \cite{Sch}) associates to $\pi$ an admissible complex representation $\phi_\pi: W'_E\to \GL(M)$ of dimension $n$. We set $L(s,\pi,\As)=L(s,\phi_\pi,\As)$, $\epsilon(s,\pi,\As,\psi')=\epsilon(s,\phi_\pi,\As,\psi')$ and $\gamma(s,\pi,\As,\psi')=\gamma(s,\phi_\pi,\As,\psi')$ and call them the Asai $L$-function, $\epsilon$-factor and $\gamma$-factor of $\pi$ respectively. We have $\det \As(\phi_\pi)=(\omega_{\pi})^n \eta_{E/F}^{n(n-1)/2}$ and therefore
\begin{align}\label{eq 3 As L functions}
\displaystyle \epsilon(s,\pi,\As,\psi'_\lambda)=\omega_{\pi}(\lambda)^n\lvert \lambda\rvert_F^{n^2(s-\frac{1}{2})} \eta_{E/F}(\lambda)^{\frac{n(n-1)}{2}}\epsilon(s,\pi,\As,\psi')
\end{align}
for every $\lambda\in F^\times$. Let $z\in \bR_+^*$. We also set
$$\displaystyle L(s,\pi,\As)=L(s,\As(\phi_\pi)).$$

Assume now that we are in the split case. Let $\pi=\pi_1\boxtimes \pi_2$ be an irreducible representation of $G_n(E)=G_n(F)\times G_n(F)$. Let $\phi_1,\phi_2:W'_F\to \GL_n(\C)$ be the admissible complex representations associated to $\pi_1$ and $\pi_2$ respectively by the local Langlands correspondence. Then, we set
$$\displaystyle L(s,\pi,\As)=L(s,\pi_1\times \pi_2)=L(s,\phi_1\otimes \phi_2),\;\; \epsilon(s,\pi,\As,\psi')=\epsilon(s,\pi_1\times \pi_2,\psi')=\epsilon(s,\phi_1\otimes \phi_2,\psi')$$
and
$$\displaystyle \gamma(s,\pi,\As,\psi')=\gamma(s,\pi_1\times \pi_2,\psi')=\gamma(s,\phi_1\otimes \phi_2,\psi').$$

We now return to the inert case. Let $P=MU$ be a standard parabolic subgroup of $G_n$ and $n_1,\ldots,n_k\in \bN^*$ so that
$$\displaystyle M(E)=G_{n_1}(E)\times\ldots\times G_{n_k}(E).$$
For each $1\leqslant i\leqslant k$, let $\sigma_i\in \Pi_2(G_{n_i}(E))$ and set $\sigma=\sigma_1\boxtimes \ldots\boxtimes \sigma_k\in \Pi_2(M(E))$. We identify $\cA_{M,\C}^*$ with $\C^k$ so that for every $(\lambda_1,\ldots,\lambda_k)\in \C^k$, we have
$$\displaystyle \sigma_\lambda=\sigma_1 \lvert \det \rvert_E^{\lambda_1}\boxtimes \ldots\boxtimes \sigma_k \lvert \det \rvert_E^{\lambda_k}.$$
We set $\pi_\lambda=i_{P(E)}^{G_n(E)}(\sigma_\lambda)$ for every $\lambda\in \cA_{M,\C}^*$. Then, $\pi_\lambda$ is irreducible for almost every $\lambda$. For such $\lambda$, $L(s,\pi_\lambda,\As)$, $\epsilon(s,\pi_\lambda,\As,\psi')$ and $\gamma(s,\pi_\lambda,\As,\psi')$ are defined as above. For the remaining $\lambda$'s, we define these factors as the one associated to the unique irreducible subquotient $\pi^0_\lambda$ of $\pi_{\lambda}$ which is the Langlands quotient of $i_{Q(E)}^{G_n(E)}(\sigma_{\lambda})$ for any parabolic $Q$ with Levi $M$ for which $\lambda$ is positive (in the large sense).

\begin{lem}\label{lem 2 As L functions}
\begin{enumerate}[(i)]
\item For every $\lambda=(\lambda_1,\ldots,\lambda_k)\in \cA_{M,\C}^*$ we have
$$\displaystyle L(s,\pi_{\lambda},\As)=\prod_{i=1}^k L(s+2\lambda_i,\sigma_i,\As)\prod_{1\leqslant i<j\leqslant k} L(s+\lambda_i+\lambda_j,\sigma_i\times \sigma_j)$$
and
$$\displaystyle \epsilon(s,\pi_{\lambda},\As,\psi')=\lambda_{E/F}(\psi')^{\frac{k(k-1)}{2}}\prod_{i=1}^k \epsilon(s+2\lambda_i,\sigma_i,\As,\psi')\prod_{1\leqslant i<j\leqslant k} \epsilon(s+\lambda_i+\lambda_j,\sigma_i\times \sigma_j,\psi'_E).$$
\item Set $U=\{(\lambda_1,\ldots,\lambda_k)\in \cA_{M,\C}^*\mid \lvert \Re(\lambda_i)\rvert<\frac{1}{4} \}$. Then, for every $s\in \C$ with $\Re(s)=\frac{1}{2}$ the function $\lambda\in U\mapsto \gamma(s,\pi_{\lambda},\As,\psi')$ is holomorphic. Moreover, for $\pi$ nearly tempered, $L(s,\pi,\As)$ has no poles in $\{\Re(s)\geqslant \frac{1}{2} \}$.
\item For any $z\in \bR_+^*$ the function $(s,\lambda)\mapsto L(s,\pi_\lambda,\As)^{-1}$ is holomorphic. Moreover, it is of finite order in vertical strips in the first variable locally uniformly in the second variable (see Section \ref{Section automatic holomorphic cont}).
\end{enumerate}
\end{lem}

\noindent\ul{Proof}:
\begin{enumerate}[(i)]
\item By compatibility of the Langlands correspondence with parabolic induction and unramified twists, for every $\lambda=(\lambda_1,\ldots,\lambda_k)\in \cA_{M,\C}^*$ the Langlands parameter of $\pi_\lambda$ (or rather $\pi^0_\lambda$) is given by
$$\displaystyle \phi_{\pi_{\lambda}}=\phi_{\sigma_1}\otimes \lvert .\rvert_E^{\lambda_1}\oplus\ldots\oplus \phi_{\sigma_k}\otimes \lvert .\rvert_E^{\lambda_k}.$$
Therefore, the identities of the lemma follow directly from \ref{eq 1 As L functions} and \ref{eq 2 As L functions}.
\item is a consequence of (i) and the fact that for $\pi_1$, $\pi_2$ tempered the $L$-functions $L(s,\pi_1,\As)$ and $L(s,\pi_1\times \pi_2)$ have no poles in $\{\Re(s)>0 \}$.
\item This is a consequence of (i) and the fact that the Gamma function decays exponentially in vertical strips.  $\blacksquare$
\end{enumerate}

\subsection{Definition and convergence of the local Zeta integrals}\label{Section defn and convergence}

Let $\pi\in \Irr_{\gen}(G_n(E))$. For every $W\in \cW(\pi,\psi_n)$, $\phi\in \cS(F^n)$ and $s\in \C$ we define, whenever convergent, a Zeta integral
$$\displaystyle Z(s,W,\phi)=\int_{N_n(F)\backslash G_n(F)} W(h)\phi(e_nh)\lvert \det h\rvert_F^s dh.$$
Recall that for every $C\in \bR$ we have set $\cH_{>C}=\{z\in \C\mid \Re(z)>C \}$ (see Section \ref{Section automatic holomorphic cont}). By Corollary \ref{cor good sections Whittaker models}, there exists $\mu\in \cA^*$ such that $\cW(\pi,\psi_n)\subseteq \cC_{\mu}(N_n(E)\backslash G_n(E),\psi_n)$. Thus, the next lemma will show that the above Zeta integrals are at least convergent in some right half-plane. 

\begin{lem}\label{lem 1 conv Zeta integrals}
Let $\mu\in \cA^*$ and $\phi\in \cS(F^n)$. Then, for every $W\in \cC_{\mu}(N_n(E)\backslash G_n(E),\psi_n)$ the integral
$$\displaystyle Z(s,W,\phi)=\int_{N_n(F)\backslash G_n(F)} W(h)\phi(e_nh)\lvert \det h\rvert_F^s dh$$
is absolutely convergent for all $s\in \cH_{>-2\min(\mu)}$. Moreover, the function $s\mapsto Z(s,W,\phi)$ is holomorphic and bounded in vertical strips. More pecisely, for every vertical strip $V\subseteq \cH_{>-2\min(\mu)}$ there exists a continuous semi-norm $p_{V,\phi}$ on $\cC_{\mu}(N_n(E)\backslash G_n(E),\psi_n)$ such that
$$\displaystyle \left\lvert Z(s,W,\phi)\right\rvert\leqslant p_{V,\phi}(W)$$
for every $W\in \cC_{\mu}(N_n(E)\backslash G_n(E),\psi_n)$ and $s\in V$.
\end{lem}

\noindent\ul{Proof}: Let $W\in \cC_{\mu}(N_n(E)\backslash G_n(E),\psi_n)$. By the Iwasawa decomposition $G_n(F)=N_n(F)A_n(F)K_n$, we need to show the convergence of
$$\displaystyle \int_{K_n}\int_{A_n(F)} \lvert W(ak)\rvert \lvert \phi(e_nak)\rvert \lvert \det a\rvert_F^t \delta_n(a)^{-1}da dk$$
for $t>-2\min(\mu)$. Since for each $R>0$ we have $\lvert \phi(e_nak)\rvert\ll (1+\lvert a_n\rvert)^{-R}$ for all $a\in A_n(F)$ and $k\in K_n$, by Lemma \ref{lem 1 HCS spaces} we are reduced to show the convergence of
\[\begin{aligned}
\displaystyle & \int_{A_n(F)} \prod_{i=1}^{n-1}(1+\lvert \frac{a_i}{a_{i+1}}\rvert)^{-R}(1+\lvert a_n\rvert)^{-R}\delta_{n,E}(a)^{1/2} a^{\lvert \mu\rvert} \lvert \det a\rvert_F^t \delta_n(a)^{-1}da= \\
 & \int_{A_n(F)} \prod_{i=1}^{n-1}(1+\lvert \frac{a_i}{a_{i+1}}\rvert)^{-R}(1+\lvert a_n\rvert)^{-R} \prod_{i=1}^n \lvert a_i\rvert^{t+2\lvert \mu\rvert_i}da
\end{aligned}\]
locally uniformly for $\frac{R}{n}-2\min(\mu)>t>-2\min(\mu)$ where $\lvert \mu\rvert_1,\ldots, \lvert \mu\rvert_n$ denote the coordinates of $\lvert \mu\rvert$. This last fact follows from the elementary inequality
$$\displaystyle \prod_{i=1}^{n-1}(1+\lvert \frac{a_i}{a_{i+1}}\rvert)^{-R}(1+\lvert a_n\rvert)^{-R}\leqslant \prod_{i=1}^n (1+\lvert a_i\rvert)^{-R/n}$$
together with the convergence of the integral
$$\displaystyle \int_{F^\times} (1+\lvert x\rvert)^{-R/n} \lvert x\rvert^{t+r}d^\times x$$
locally uniformly for $\frac{R}{n}-r>t>-r$ for every $r\in \bR$. $\blacksquare$

We recall that we have introduced in Section \ref{section notation} the notion of nearly tempered representation of $G_n(E)$.

\begin{lem}\label{lem 2 conv Zeta integrals}
Assume that $\pi$ is a generic irreducible representation of $G_n(E)$ which is nearly tempered. Then, there exists $\epsilon>0$ such that for every $W\in \cW(\pi,\psi_n)$ and $\phi\in \cS(F^n)$ the integral defining $Z(s,W,\phi)$ converges absolutely on $\cH_{\frac{1}{2}-\epsilon}$ and defines a holomorphic function bounded in vertical strips there.
\end{lem}

\noindent\ul{Proof}: Let $W\in \cW(\pi,\psi_n)$ and $\phi\in \cS(F^n)$. By assumption, there exists a parabolic subgroup $P=MU$ of $G_n$, a discrete series $\sigma$ of $M(E)$ and $\lambda=(\lambda_1,\ldots,\lambda_n)\in \cA_{M,\C}^*\subseteq \cA_{\C}^*=\C^n$ satisfying $\lvert \Re(\lambda_i)\rvert<\frac{1}{4}$ for every $1\leqslant i\leqslant n$ such that $\pi\simeq i_{P(E)}^{G_n(E)}(\sigma_{\lambda})$. Let $\rho=(\frac{n-1}{2},\ldots,\frac{1-n}{2})\in \cA^*$ be half the sum of the roots of $A_n$ in $B_n$. Then, for every $\eta>0$ we have $\lvert \Re(\lambda)\rvert\prec \lvert \Re(\lambda)\rvert+\eta \rho$. Therefore, by Proposition \ref{prop 1 uniform bounds Whittaker functions}(i) (and standard properties of the Jacquet's functional \cite[Sect. 15.4]{Wall2}) we have
$$\displaystyle \cW(\pi,\psi_n)\subseteq \cC_{\lvert \Re(\lambda)\rvert+\eta \rho}(N_n(E)\backslash G_n(E),\psi_n)$$
for any $\eta>0$. By the previous lemma it follows that for any $\eta>0$, $Z(s,W,\phi)$ converges absolutely on $\cH_{-2\min(\lvert \Re(\lambda)\rvert+\eta \rho)}$ and defines a holomorphic function there. Since
$$\displaystyle \cH_{-2\min(\lvert \Re(\lambda)\rvert)}=\bigcup_{\eta>0} \cH_{-2\min(\lvert \Re(\lambda)\rvert+\eta \rho)}$$
we deduce that $Z(s,W,\phi)$ converges absolutely on $\cH_{-2\min(\lvert \Re(\lambda)\rvert)}$ and defines a holomorphic function there. Finally, the inequalities satisfied by $\lambda$ imply that $\cH_{-2\min(\lvert \Re(\lambda)\rvert)}=\cH_{\frac{1}{2}-\epsilon}$ for some $\epsilon>0$ and the lemma follows. $\blacksquare$

\vspace{2mm}

We end this section with the following non-vanishing result.

\begin{lem}\label{lem nonzero Zeta}
For every $s_0\in \C$, there exist finite families $W_i\in \cW(\pi,\psi_n)$ and $\phi_i\in \cS(F^n)$ indexed by $i\in I$ such that the function $s\mapsto \sum_{i\in I}Z(s,W_i,\phi_i)$ (which is defined on some right-half plane) admits a holomorphic continuation to $\C$ which is non-vanishing at $s_0$.
\end{lem}

\noindent\ul{Proof}: Let $P_n$ be the mirabolic subgroup of $G_n$ (i.e. the subgroup of elements $g\in G_n$ with last row $(0,\ldots,0,1)$), $U_n$ be the unipotent radical of $P_n$ and $\overline{U}_n={}^tU_n$. Then, we have the decomposition $G_n(F)=P_n(F)Z_n(F)\overline{U}_n(F)$ and correspondingly for Haar measures $dh=\lvert \det p\rvert_F^{-1}d_r pdz d\overline{u}$ where $d_r p$ denotes a right Haar measure on $P_n(F)$. Thus, by Lemma \ref{lem 1 conv Zeta integrals}, for $\Re(s)\gg 1$ the expression defining $Z(s,W,\phi)$ is absolutely convergent and we have
\[\begin{aligned}
\displaystyle Z(s,W,\phi) & =\int_{Z_n(F)\times \overline{U}_n(F)} \int_{N_n(F)\backslash P_n(F)} W(pz\overline{u}) \lvert \det p\rvert_F^{s-1} d_rp \phi(e_nz\overline{u}) \lvert \det z\rvert^s_F dzd\overline{u} \\
 & =\int_{Z_n(F)\times \overline{U}_n(F)} \int_{N_n(F)\backslash P_n(F)} W(p\overline{u}) \lvert \det p\rvert_F^{s-1} d_rp \phi(e_nz\overline{u}) \omega_{\pi}(z)\lvert \det z\rvert^s_F dzd\overline{u}
\end{aligned}\]
for every $W\in \cW(\pi,\psi_n)$ and $\phi\in \cS(F^n)$. Let $\varphi_Z\in C_c^\infty(Z_n(F))$ and $\varphi_{\overline{U}}\in C_c^\infty(\overline{U}_n(F))$. Then, there exists a unique $\phi=\phi_{\varphi_Z,\varphi_{\overline{U}}}\in C_c^\infty(F^n)$ such that $\phi(e_nz\overline{u})=\varphi_Z(z)\varphi_{\overline{U}}(\overline{u})$ for every $(z,u)\in Z_n(F)\times \overline{U}_n(F)$. For such a $\phi$, the above identity becomes
$$\displaystyle Z(s,W,\phi)=\int_{N_n(F)\backslash P_n(F)} \left(R(\varphi_{\overline{U}})W\right)(p) \lvert \det p\rvert_F^{s-1} d_rp \int_{Z_n(F)} \varphi_Z(z) \omega_{\pi}(z) \lvert \det z\rvert^s dz,\;\;\; \Re(s)\gg 1$$
for every $W\in \cW(\pi,\psi_n)$. Let $f\in C_c^\infty(N_n(E)\backslash P_n(E),\psi_n)$. By \cite{GK}, \cite[Proposition 5]{Jac3} and \cite{Kem}, there exists $W_0\in \cW(\pi,\psi_n)$ whose restriction to $P_n(E)$ coincides with $f$. By Dixmier-Malliavin \cite{DM} (in the Archimedean case), there exist finite families $(W_i)_{i\in I}$ and $(\varphi_{\overline{U},i})_{i\in I}$ of elements in $\cW(\pi,\psi_n)$ and $C_c^\infty(\overline{U}_n(F))$ respectively such that
$$\displaystyle W_0=\sum_{i\in I} R(\varphi_{\overline{U},i})W_i.$$
Choose $\varphi_Z\in C_c^\infty(Z_n(F))$ such that
$$\displaystyle \int_{Z_n(F)} \varphi_Z(z) \omega_{\pi}(z) \lvert \det z\rvert^{s_0} dz\neq 0$$
(Notice that the above integral is absolutely convergent for any complex value of $s_0$). Set $\phi_i=\phi_{\varphi_Z, \varphi_{\overline{U},i}}$ for every $i\in I$. Then, by the above, for $\Re(s)\gg 1$ we have
\[\begin{aligned}
\displaystyle \sum_{i\in I} Z(s,W_i,\phi_i) & =\sum_{i\in I} \int_{N_n(F)\backslash P_n(F)} \left(R(\varphi_{\overline{U},i})W_i\right)(p) \lvert \det p\rvert_F^{s-1} d_rp \int_{Z_n(F)} \varphi_Z(z) \omega_{\pi}(z) \lvert \det z\rvert^s dz \\
 & =\int_{N_n(F)\backslash P_n(F)} W_0(p) \lvert \det p\rvert_F^{s-1} d_rp \int_{Z_n(F)} \varphi_Z(z) \omega_{\pi}(z) \lvert \det z\rvert^s dz \\
 & =\int_{N_n(F)\backslash P_n(F)} f(p) \lvert \det p\rvert_F^{s-1} d_rp \int_{Z_n(F)} \varphi_Z(z) \omega_{\pi}(z) \lvert \det z\rvert^s dz.
\end{aligned}\]
The above integrals are convergent for any $s\in \C$ uniformly on compacta and therefore the resulting expression defines a holomorphic function on $\C$. Moreover, we can certainly choose $f$ so that
$$\displaystyle \int_{N_n(F)\backslash P_n(F)} f(p) \lvert \det p\rvert_F^{s_0-1} d_rp\neq 0.$$
By our choice of $\varphi_Z$ this implies that the holomorphic continuation of $\sum_{i\in I} Z(s,W_i,\phi_i)$ does not vanish at $s_0$. $\blacksquare$

\subsection{Local functional equation: the split case}\label{Section split case}

In this section we assume that we are in the split case i.e. $E=F\times F$. Let $\pi\in \Irr_{\gen}(G_n(E))$. Then $\pi=\pi_1\boxtimes \pi_2$ for some $\pi_1,\pi_2\in \Irr(G_n(F))$. In the case where $\psi(x,y)=\psi'(x)\psi'(-y)$ for $(x,y)\in E$ (i.e. $\tau=(1,-1)$) and $W=W_1\otimes W_2$ for some $W_1\in \cW(\pi_1,\psi'_n)$, $W_2\in \cW(\pi_2,{\psi'_n}^{-1})$, for $\phi\in \cS(F^n)$ the Zeta integral $Z(s,W,\phi)$ belongs to a family of expressions studied by Jacquet-Piatetskii-Shapiro-Shalika and Jacquet in \cite{JPSS} and \cite{Jac}. By the main results of those references together with some of the characterizing properties of the local Langlands correspondence for $\GL_n$ (\cite{HT}, \cite{Hen}, \cite{Sch}), in this situation $Z(s,W,\phi)$ admits a meromorphic continuation to $\C$ satisfying the functional equation
\[\begin{aligned}
\displaystyle  Z(1-s,\widetilde{W},\widehat{\phi}) & =\omega_{\pi_2}(-1)^{n-1}\gamma(s,\pi_1\times \pi_2,\psi')Z(s,W,\phi) \\
 & =\omega_\pi(\tau)^{n-1}\lvert \tau\rvert_E^{\frac{n(n-1)}{2}(s-1/2)}\lambda_{E/F}(\psi')^{-\frac{n(n-1)}{2}}\gamma(s,\pi,\As,\psi')Z(s,W,\phi).
\end{aligned}\]
Let us remark here that our conventions are slightly different to the ones in \cite{Jac}: the Fourier transform $\phi\mapsto \widehat{\phi}$ in {\em loc. cit.} is normalized using the character $\psi'^{-1}$ rather than $\psi'$ (see \cite[\S 2]{Jac}) which results in composing the one used in this paper with the involution $x\in F^n\mapsto -x$. On the other hand, in the functional equation of \cite[Theorem 2.1]{Jac} the term $\omega_{\pi_2}(-1)^{n-1}$ is replaced by $\omega_{\pi_1}(-1)^{n-1}$. Finally, in \cite{Jac} the local $\gamma$-factors are also normalized differently: although in this paper we have followed the normalization of \cite{Ta} and \cite{GR}, in {\em loc. cit.} the convention is opposite and what we denote here by $\gamma(s,\pi_1\times \pi_2,\psi')$ corresponds to the factor $\gamma(s,\pi_1\times \pi_2,\psi'^{-1})$ in \cite{Jac} (see in particular the appendix \cite[Section 16]{Jac}). Taking all of these into account, it is a straightforward exercise to see that the above functional equation is equivalent to the one given in \cite[Theorem 2.1]{Jac}.

In the $p$-adic case the meromorphic continuation of $s\mapsto Z(s,W,\phi)$ and the above equality extend to any $W\in \cW(\pi,\psi)$, by linearity, whereas in the Archimedean case such an extension follows from the remark after Theorem 2.3 of \cite{Jac}. We record this as a theorem by removing the assumption on $\psi$.

\begin{theo}\label{theo 1 split}
Assume that $E=F\times F$. Let $\pi\in \Irr_{\gen}(G_n(E))$. Then, for every $W\in \cW(\pi,\psi_n)$ and $\phi\in \cS(F^n)$ the function $s\mapsto Z(s,W,\phi)$ has a meromorphic extension to $\C$ satisfying the functional equation
\begin{align}\label{eq 0 split}
\displaystyle  Z(1-s,\widetilde{W},\widehat{\phi})=\omega_\pi(\tau)^{n-1}\lvert \tau\rvert_E^{\frac{n(n-1)}{2}(s-1/2)}\lambda_{E/F}(\psi')^{-\frac{n(n-1)}{2}}\gamma(s,\pi,\As,\psi')Z(s,W,\phi).
\end{align}
\end{theo}

\noindent\ul{Proof}: By the above discussion, the theorem holds when $\psi(x,y)=\psi'(x)\psi'(-y)$ for every $(x,y)\in E$. We therefore just need to study the effect of replacing $\psi$ by $\psi_\lambda$ for some $\lambda\in F^\times$ where $\psi_\lambda(z)=\psi(\lambda z)$ for every $z\in E$. Doing so amounts to replacing $\tau$ by $\lambda \tau$. Define the generic character $\psi_{n,\lambda}$ as $\psi_n$ using $\psi_\lambda$ instead of $\psi$. Then, there is an isomorphism $\cW(\pi,\psi_n)\to \cW(\pi,\psi_{n,\lambda})$, $W\mapsto W_\lambda$ given by
$$\displaystyle W_\lambda(g)=W(a(\lambda)g),\;\;\; g\in G_n(E)$$
where $a(\lambda)=\begin{pmatrix} \lambda^{n-1} & & \\ & \ddots & \\ & & 1 \end{pmatrix}$. By the change of variable $h\mapsto a(\lambda)^{-1}h$ we have
\begin{align}\label{eq 1 split}
\displaystyle Z(s,W_\lambda,\phi)=\delta_n(a(\lambda))\lvert \det a(\lambda)\rvert_F^{-s} Z(s,W,\phi),\;\;\; \Re(s)\gg 1.
\end{align}
This already shows that $Z(s,W_\lambda,\phi)$ has a meromorphic continuation if and only if $Z(s,W,\phi)$ has one. On the other hand, for all $g\in G_n(E)$ we have
\[\begin{aligned}
\displaystyle \widetilde{(W_\lambda)}(g)=W(a(\lambda)w_n{}^t g^{-1})=W(w_n{}^t(w_na(\lambda)^{-1}w^{-1}_ng)^{-1})=\widetilde{W}(w_na(\lambda)^{-1}w^{-1}_ng)=\omega_{\pi}(\lambda)^{n-1}(\widetilde{W})_\lambda(g)
\end{aligned}\]
so that $\widetilde{W_\lambda}=\omega_{\pi}(\lambda)^{n-1}(\widetilde{W})_\lambda$ and finally (by \ref{eq 1 split})
\begin{align}\label{eq 2 split}
\displaystyle Z(1-s,\widetilde{W_\lambda},\widehat{\phi})=\omega_{\pi}(\lambda)^{n-1}\delta_n(a(\lambda))\lvert \det a(\lambda)\rvert_F^{s-1} Z(1-s,\widetilde{W},\widehat{\phi}).
\end{align}
From \ref{eq 1 split} and \ref{eq 2 split} it follows that if $Z(s,W,\phi)$ and $Z(s,\widetilde{W},\widehat{\phi})$ satisfy \ref{eq 0 split} then
\[\begin{aligned}
\displaystyle Z(1-s,\widetilde{W_\lambda},\widehat{\phi}) & =\omega_{\pi}(\lambda)^{n-1}\lvert \det a(\lambda)\rvert_F^{2s-1}\omega_\pi(\tau)^{n-1}\lvert \tau\rvert_E^{\frac{n(n-1)}{2}(s-1/2)}\lambda_{E/F}(\psi')^{-\frac{n(n-1)}{2}}\gamma(s,\pi_1\times \pi_2,\psi') \\
 & \times Z(s,W_\lambda,\phi) \\
& =\omega_\pi(\lambda \tau)^{n-1}\lvert \lambda\tau\rvert_E^{\frac{n(n-1)}{2}(s-1/2)}\lambda_{E/F}(\psi')^{-\frac{n(n-1)}{2}}\gamma(s,\pi_1\times \pi_2,\psi')Z(s,W_\lambda,\phi)
\end{aligned}\]
which is precisely the functional equation for $\psi$ replaced by $\psi_\lambda$. $\blacksquare$

\subsection{Local functional equation: the inert case}\label{Section inert case}

In this Section we assume that $E/F$ is a quadratic field extension. We are going to state two theorems which are the main results of this paper. These theorems will be proved in Sections \ref{appendix proof of the theorem} and \ref{Section end of proof of main theorems}. The first result is the exact analog of Theorem \ref{theo 1 split} in the inert case:

\begin{theo}\label{theo 1 inert}
Assume that $E/F$ is a quadratic field extension. Let $\pi\in \Irr_{\gen}(G_n(E))$. Then, for every $W\in \cW(\pi,\psi_n)$ and $\phi\in \cS(F^n)$ the function $s\mapsto Z(s,W,\phi)$ has a meromorphic extension to $\C$ and satisfies the functional equation
\begin{align}\label{eq 0 inert}
\displaystyle  Z(1-s,\widetilde{W},\widehat{\phi})=\omega_\pi(\tau)^{n-1}\lvert \tau\rvert_E^{\frac{n(n-1)}{2}(s-1/2)}\lambda_{E/F}(\psi')^{-\frac{n(n-1)}{2}}\gamma(s,\pi,\As,\psi')Z(s,W,\phi).
\end{align}
\end{theo}

\begin{theo}\label{theo 2 inert}
Assume that $E/F$ is a quadratic field extension. Let $\pi\in \Irr_{\gen}(G_n(E))$. Then, for every $W\in \cW(\pi,\psi_n)$ and $\phi\in \cS(F^n)$ the function
$$\displaystyle s\mapsto \frac{Z(s,W,\phi)}{L(s,\pi,\As)}$$
is holomorphic and of finite order in vertical strips. Moreover, if $\pi$ is nearly tempered for every $s_0\in \C$ there exist $W\in \cW(\pi,\psi_n)$ and $\phi\in \cS(F^n)$ for which the above function does not vanish at $s_0$.
\end{theo}

As in the introduction, we recall here that a big part of Theorems \ref{theo 1 inert} and \ref{theo 2 inert} was already known in the $p$-adic case. More precisely, in this case the functional equation of Theorem \ref{theo 1 inert} was established by Flicker \cite[Appendix]{Fli} and Kable \cite[Theorem 3]{Kab} with a possibly different $\gamma$-factor $\gamma^{\RS}(s,\pi,\As,\psi')=\epsilon^{\RS}(s,\pi,\As,\psi')\frac{L^{\RS}(1-s,\pi^\vee,\As)}{L^{\RS}(s,\pi,\As)}$. Here, by its very definition, $L^{\RS}(s,\pi,\As)$ is the inverse of a polynomial in $q^{-s}$ with $1$ as constant term such that $\frac{Z(s,W,\phi)}{L^{\RS}(s,\pi,\As)}$ is always holomorphic and non-vanishing at any given point $s_0\in \C$ for suitable choice of the pair $(W,\phi)\in \cW(\pi,\psi_n)\times \cS(F^n)$. The equality $L^{\RS}(s,\pi,\As)=L(s,\pi,\As)$ was then proved first by Anandavardhanan-Rajan \cite{AR} for $\pi$ square-integrable and then by Matringe \cite{Mat} in general. On the other hand, the identity of $\epsilon$-factors $\epsilon^{\RS}(s,\pi,\As,\psi')=\epsilon(s,\pi,\As,\psi')$ was previously only known when $\pi$ is supercuspidal (\cite{AKMSS}) or $n=2$ (\cite{CCI}). Therefore, this identity of $\epsilon$-factors is the only truly new result when $E$ and $F$ are $p$-adic. However, we won't rely on those previous works in the proof, except for the global functional equation of Flicker and Kable (see Section \ref{Global functional equation}), and moreover we will be able to treat mostly uniformly non-Archimedean and Archimedean fields.

The aim of the next lemma is to check that Theorems \ref{theo 1 inert} and \ref{theo 2 inert} do not depend on the choices of $\psi$ and $\psi'$.

\begin{lem}\label{lem 1 inert}
If Theorems \ref{theo 1 inert} and \ref{theo 2 inert} hold for one pair $(\psi,\psi')$ of nontrivial additive characters of $E$, $F$ respectively with $\psi$ trivial on $F$ then they hold for any such pair.
\end{lem}

\noindent\ul{Proof}: The independence on the choice of $\psi$ can be proved exactly the same way as in the proof of Theorem \ref{theo 1 split}. Clearly Theorem \ref{theo 2 inert} and the meromorphic extension of $s\mapsto Z(s,W,\phi)$ for every $W\in \cW(\pi,\psi_n)$ and $\phi\in \cS(F^n)$ do not depend on the choice of $\psi'$. Thus, it only remains to show the independence of the functional equation \ref{eq 0 inert} on the choice of $\psi'$. If we replace $\psi'$ by $\psi'_\lambda$ defined by $\psi'_\lambda(x)=\psi'(\lambda x)$ for every $x\in F$ for some $\lambda\in F^\times$ but keep $\psi$ fixed then we have to replace $\tau$ by $\lambda^{-1}\tau$. Therefore, the functional equation \ref{eq 0 inert} for $\psi'$ replaced by $\psi'_\lambda$ reads
\begin{align}\label{eq 1 inert}
\displaystyle  Z(1-s,\widetilde{W},\widehat{\phi}^{\psi'_\lambda})=\omega_\pi(\lambda^{-1}\tau)^{n-1}\lvert \lambda^{-1}\tau\rvert_E^{\frac{n(n-1)}{2}(s-1/2)}\lambda_{E/F}(\psi'_\lambda)^{-\frac{n(n-1)}{2}}\gamma(s,\pi,\As,\psi'_\lambda)Z(s,W,\phi)
\end{align}
for every $W\in \cW(\pi,\psi_n)$ and $\phi\in \cS(F^n)$ where $\widehat{\phi}^{\psi'_\lambda}$ stands for the Fourier transform of $\phi$ with respect to $\psi'_\lambda$ (rather than $\psi'$) and the corresponding autodual Haar measure on $F^n$. We have the relation $\widehat{\phi}^{\psi'_\lambda}(v)=\lvert \lambda\rvert_F^{n/2}\widehat{\phi}(\lambda v)$ ($v\in F^n$). Therefore, by the change of variable $h\mapsto \lambda^{-1}h$, we have
\begin{align}\label{eq 2 inert}
\displaystyle Z(1-s,\widetilde{W},\widehat{\phi}^{\psi'_\lambda})=\lvert \lambda\rvert_F^{n(s-1/2)} \omega_{\pi}(\lambda)Z(1-s,\widetilde{W},\widehat{\phi}).
\end{align}
On the other hand,
\begin{align}\label{eq 3 inert}
\displaystyle \lambda_{E/F}(\psi'_\lambda)=\eta_{E/F}(\lambda)\lambda_{E/F}(\psi')
\end{align}
whereas by \ref{eq 3 As L functions},
\begin{align}\label{eq 4 inert}
\displaystyle \gamma(s,\pi,\As,\psi'_\lambda)=\omega_{\pi}(\lambda)^n \lvert \lambda\rvert_F^{n^2(s-1/2)} \eta_{E/F}(\lambda)^{n(n-1)/2}\gamma(s,\pi,\As,\psi').
\end{align}
Combining \ref{eq 2 inert}, \ref{eq 3 inert} and \ref{eq 4 inert} we see that the functional equation \ref{eq 1 inert} for $(\psi,\psi'_\lambda)$ reduces to the one for $(\psi,\psi')$ (that is \ref{eq 0 inert}). $\blacksquare$

The next lemma is straightforward and allows to restrict to representations with unitary central character.

\begin{lem}\label{lem unramified twists}
Let $\pi$ be a generic irreducible representation of $G_n(E)$ and $\lambda\in \cA_{G_n,\C}^*$. Then if Theorems \ref{theo 1 inert} and \ref{theo 2 inert} hold for $\pi_{\lambda}$ they also hold for $\pi$.
\end{lem}

\subsection {Unramified computation}\label{Section unr computation}

In this Section, we consider the case where $F$ is non-Archimedean and the extension $E/F$ is either inert or split. Recall that an irreducible representation $\pi$ of $G_n(E)$ is said to be {\em unramified} if it admits a nonzero $G_n(\cO_E)$-fixed vector in which case $\pi^{G_n(\cO_E)}$ is a line. We also say that the characters $\psi'$ and $\psi$ are {\em unramified} if the maximal fractional ideals on which they are trivial are $\cO_F$ and $\cO_E$ respectively. If $\pi$ and $\psi$ are unramified there exists a unique $W\in \cW(\pi,\psi_n)^{G_n(\cO_E)}$ such that $W(1)=1$. The following unramified computation is standard and already in the literature in all but one case. We provide a proof in this missing case.

\begin{lem}\label{lem 1 unramified}
Assume that $F$ is non-Archimedean and that $\pi\in \Irr_{\gen}(G_n(E))$, $\psi'$ and $\psi$ are all unramified. Let $W\in \cW(\pi,\psi_n)^{G_n(\cO_E)}$ be normalized by $W(1)=1$ and $\phi\in \cS(F^n)$ be the characteristic function of $\cO_F^n$. Then, for $\Re(s)\gg 1$ we have
$$\displaystyle Z(s,W,\phi)=\vol(N_n(\cO_F)\backslash G_n(\cO_F))L(s,\pi,\As).$$
\end{lem}

\noindent\ul{Proof}: In the split case this is \cite[Proposition 2.3]{JS} and in the inert case when the extension $E/F$ is unramified this is \cite[Proposition 3]{Fli2}. Thus, it only remains to deal with the case where $E/F$ is a ramified field extension. Let $\varpi_F$ be an uniformizer of $F$. For any $n$-uple $\lambda=(\lambda_1,\ldots,\lambda_n)$ of integers set 
$$\displaystyle a(\lambda)=\begin{pmatrix} \varpi_F^{\lambda_1} \\ & \ddots \\ & & \varpi_F^{\lambda_n}\end{pmatrix}.$$
Then, by the Iwasawa decomposition we have (for $\Re(s)\gg 1$)
\begin{align}\label{eq 0 unramified}
\displaystyle Z(s,W,\phi)=\vol(N_n(\cO_F)\backslash G_n(\cO_F))\sum_{\lambda\in\bZ^n}W(a(\lambda))\phi(\varpi_F^{\lambda_n}e_n)\lvert \det a(\lambda) \rvert_F^s \delta_n(a(\lambda))^{-1}
\end{align}
For every $n$-uple $\lambda=(\lambda_1,\ldots,\lambda_n)$ of decreasing integers let $s_\lambda$ be the Schur function as defined in \cite[\S 3]{Fli2}. Let $t=(t_1,\ldots,t_n)\in (\C^\times)^n/\fS_n$ be the Satake parameter of $\pi$. Then, by \cite{CS}, \cite{Shin}, $W(a(\lambda))$ is zero unless $\lambda_1\geqslant \ldots \geqslant \lambda_n$ in which case it equals $\delta_{n,E}(a(\lambda))^{1/2}s_{2\lambda}(t)$. Moreover as $\phi=\mathbf{1}_{\cO_F^n}$ we have $\phi(\varpi_F^{\lambda_n}e_n)=0$ if $\lambda_n<0$ and $1$ otherwise. Therefore, by \ref{eq 0 unramified} and since $\delta_{n,E}=\delta_{n}^2$ on $A_n(F)$ we get
\[\begin{aligned}
\displaystyle \vol(N_n(\cO_F)\backslash G_n(\cO_F))^{-1}Z(s,W,\phi) & =\sum_{\substack{\lambda\in\bZ^n \\ \lambda_1\geqslant \ldots\geqslant \lambda_n\geqslant 0}}s_{2\lambda}(t)\lvert \det a(\lambda) \rvert_F^s \\
 & =\sum_{\substack{\lambda\in\bZ^n \\ \lambda_1\geqslant \ldots\geqslant \lambda_n\geqslant 0}}s_{2\lambda}(t)q_F^{-s(\lambda_1+\ldots+\lambda_n)} \\
 & =\sum_{\substack{\lambda\in\bZ^n \\ \lambda_1\geqslant \ldots\geqslant \lambda_n\geqslant 0}}s_{2\lambda}(q_F^{-s/2}t_1,\ldots,q_F^{-s/2}t_n)
\end{aligned}\]
But by \cite{M} equality 5.(a) p.77 we have
$$\displaystyle \sum_{\substack{\lambda\in\bZ^n \\ \lambda_1\geqslant \ldots\geqslant \lambda_n\geqslant 0}}s_{2\lambda}(q_F^{-s/2}t_1,\ldots,q_F^{-s/2}t_n)=\prod_{1\leqslant i\leqslant n}(1-q_F^{-s}t_i^2)^{-1}\prod_{1\leqslant i<j\leqslant n} (1-q_F^{-s}t_it_j)^{-1}.$$
On the other hand, since $\pi$ is unramified it is of the form $\chi_1\times\ldots\times \chi_n$ where the $\chi_i$ are unramified character of $E^\times$. Thus, by Lemma \ref{lem 2 As L functions}(i), we have
\[\begin{aligned}
\displaystyle L(s,\pi,\As) & =\prod_{i=1}^n L(s,\chi_i,\As)\prod_{1\leqslant i<j\leqslant n} L(s,\chi_i\chi_j)=\prod_{i=1}^n L(s,\chi_{i\mid F^\times})\prod_{1\leqslant i<j\leqslant n} L(s,\chi_i\chi_j) \\
 & =\prod_{1\leqslant i\leqslant n}(1-q_F^{-s}t_i^2)^{-1}\prod_{1\leqslant i<j\leqslant n} (1-q_F^{-s}t_it_j)^{-1}
\end{aligned}\]
and this ends the proof of the lemma. $\blacksquare$

\begin{cor}\label{cor unramified}
Assume that $E/F$ is a quadratic extension of $p$-adic fields and that $\pi\in \Irr_{\gen}(G_n(E))$ is an unramified generic representation of $G_n(E)$. Then, there exist $W\in \cW(\pi,\psi_n)$ and $\phi\in \cS(F^n)$ such that $s\mapsto Z(s,W,\phi)$ and $s\mapsto Z(s,\widetilde{W},\widehat{\phi})$ are not identically zero and admit meromorphic continuation to $\C$ satisfying
$$\displaystyle  Z(1-s,\widetilde{W},\widehat{\phi})=\omega_\pi(\tau)^{n-1}\lvert \tau\rvert_E^{\frac{n(n-1)}{2}(s-1/2)}\lambda_{E/F}(\psi')^{-\frac{n(n-1)}{2}}\gamma(s,\pi,\As,\psi')Z(s,W,\phi).$$
\end{cor}

\noindent\ul{Proof}: By Lemma \ref{lem 1 inert} (or rather its proof) the statement of the corollary does not depend on the choice of $(\psi,\psi')$. Therefore, we may assume that both $\psi'$ and $\psi$ are unramified. Then, taking for $W$ the unique element of $\cW(\pi,\psi_n)^{G_n(\cO_E)}$ normalized by $W(1)=1$ and for $\phi\in \cS(F^n)$ the characteristic function of $\cO_F^n$, by the previous lemma we have
$$\displaystyle Z(s,W,\phi)=\vol(N_n(\cO_F)\backslash G_n(\cO_F))L(s,\pi,\As),$$
$$\displaystyle Z(s,\widetilde{W},\widehat{\phi})=\vol(N_n(\cO_F)\backslash G_n(\cO_F))L(s,\widetilde{\pi},\As)$$
for $\Re(s)\gg 1$. This shows that both $s\mapsto Z(s,W,\phi)$ and $s\mapsto Z(s,\widetilde{W},\widehat{\phi})$ are not identically zero and admit meromorphic continuation to $\C$ satisfying
$$\displaystyle Z(1-s,\widetilde{W},\widehat{\phi})=\frac{L(1-s,\widetilde{\pi},\As)}{L(s,\pi,\As)}Z(s,W,\phi).$$
Thus, it only remains to check that
$$\displaystyle \epsilon(s,\pi,\As,\psi')=\omega_\pi(\tau)^{1-n}\lvert \tau\rvert_E^{\frac{n(n-1)}{2}(1/2-s)}\lambda_{E/F}(\psi')^{\frac{n(n-1)}{2}}.$$
Since $\pi$ is unramified, it is an induced representation of the form $\chi_1\times\ldots\times \chi_n$ where the $\chi_i$'s are unramified characters of $E^\times$. By Lemma \ref{lem 2 As L functions}(i) we have
\begin{align}\label{eq 2 epsilon}
\displaystyle \epsilon(s,\pi,\As,\psi')=\lambda_{E/F}(\psi')^{\frac{n(n-1)}{2}}\prod_{1\leqslant i<j\leqslant n}\epsilon(s,\chi_i\chi_j,\psi'_E)\times \prod_{i=1}^n \epsilon(s,\chi_{i\mid F^\times},\psi')
\end{align}
Since $\chi_{i\mid F^\times}$ and $\psi'$ are unramified we have $\epsilon(s,\chi_{i\mid F^\times},\psi')=1$ for all $1\leqslant i\leqslant n$. On the other hand, since $\chi_i$, $\chi_j$ and $\psi=\psi'_{E,\tau}$ are unramified we have
$$\displaystyle \epsilon(s,\chi_i\chi_j,\psi'_E)=\lvert \tau\rvert_E^{1/2-s}\chi_i(\tau)^{-1}\chi_j(\tau)^{-1}\epsilon(s,\chi_i\chi_j,\psi)=\lvert \tau\rvert_E^{1/2-s}\chi_i(\tau)^{-1}\chi_j(\tau)^{-1}$$
for every $1\leqslant i<j\leqslant n$. Therefore \ref{eq 2 epsilon} becomes
\[\begin{aligned}
\displaystyle \epsilon(s,\pi,\As,\psi') & =\lambda_{E/F}(\psi')^{\frac{n(n-1)}{2}}\lvert \tau\rvert_E^{\frac{n(n-1)}{2}(1/2-s)}\prod_{1\leqslant i<j\leqslant n}\chi_i(\tau)^{-1}\chi_j(\tau)^{-1} \\
 & =\lambda_{E/F}(\psi')^{\frac{n(n-1)}{2}}\lvert \tau\rvert_E^{\frac{n(n-1)}{2}(1/2-s)} \omega_{\pi}(\tau)^{1-n}
\end{aligned}\]
which is exactly what we wanted. $\blacksquare$

\subsection{Global Zeta integrals and their functional equation}\label{Global functional equation}

In this section we let $k'/k$ be a quadratic extension of number fields. For every place $v$ of $k$ we denote by $k_v$ the corresponding completion of $k$ and set $k'_v=k_v\otimes_k k'$. Also if $v$ is non-Archimedean, we let $\cO_v$, $\cO_{k',v}$ be the ring of integers of $k_v$ and $k'_v$ respectively. Let $\bA=\prod'_{v} k_v$ and $\bA_{k'}=\bA\otimes_k k'=\prod_v' k'_v$ be the ad\`ele rings of $k$ and $k'$ respectively and $\lvert .\rvert$ be the normalized absolute value on $\bA$. We will also denote by $\bA_{k',f}=\prod'_{v\not\mid \infty} k'_v$ the finite ad\`eles of $k'$ and by $k'_\infty=\prod_{v\mid \infty} k'_v$ the product of the Archimedean completions of $k'$ (so that $\bA_{k'}=\bA_{k',f}\times k'_\infty$). Let $\Psi'$ and $\Psi$ be nontrivial additive characters of $k\backslash \bA$ and $k'\backslash \bA_{k'}$ respectively with $\Psi$ {\em trivial} on $k\backslash \bA$. For every place $v$ we let $\Psi_{v}$ be the local component of $\Psi$ at $v$ to which we associate a generic character $\Psi_{n,v}:N_n(k'_v)\to \mathbb{S}^1$ as before. Then $\Psi_n=\prod_v \Psi_{n,v}$ define a character of $N_n(\bA_{k'})$ which is at the same time trivial on $N_n(k')$ and $N_n(\bA)$. Let $Z_\infty$ be the connected component (for the usual topology) of $Z_n(\R)$. By the diagonal embedding $\R\hookrightarrow k'_\infty$, we consider $Z_\infty$ as a central subgroup of $G_n(k'_\infty)$. For $\omega$ a continuous unitary character of $Z_\infty$ we denote by $\cA_{\cusp}(Z_\infty G_n(k')\backslash G_n(\bA_{k'}),\omega)$ the space of smooth functions $\varphi: G_n(k')\backslash G_n(\bA_{k'})\to \C$ (here smooth means that $\varphi$ is right invariant by an open subgroup of $G_n(\bA_{k',f})$ and $C^\infty$ in the $G(k'_\infty)$ component) with all their derivatives of moderate growth in the sense of \cite[\S I.2.3]{MW}, having central character $\omega$ (i.e. $\varphi(zg)=\omega(z)\varphi(g)$ for every $(z,g)\in Z_\infty\times G_n(\bA_{k'})$) and satisfying
$$\displaystyle \int_{N(k')\backslash N(\bA_{k'})}\varphi(ug)du=0,\;\;\; g\in G_n(\bA_{k'})$$
for every parabolic subgroup $P=MN$ of $G_n$. By a cuspidal automorphic representation of $G_n(\bA_{k'})$ we mean a closed topologically irreducible subrepresentation of the space $\cA_{\cusp}(Z_\infty G_n(k')\backslash G_n(\bA_{k'}),\omega)$ for some character $\omega$. If $\Pi$ is such a cuspidal automorphic representation of $G_n(\bA_{k'})$ then there exists an isomorphism
\begin{align}\label{eq 1 globalization}
\displaystyle \Pi\simeq \widehat{\bigotimes_{v\mid \infty}}\Pi_v\otimes \bigotimes'_{v\not\mid \infty}\Pi_v
\end{align}
where for every place $v$ of $k$, $\Pi_v$ is an irreducible representation of $G_n(k'_v)$ in the sense of Section \ref{Section representations}. Moreover, the representation $\Pi_v$ is generic for every place $v$, unramified for almost all places $v$ and the restricted tensor product above is taken with respect to a family of (almost everywhere defined) $G_n(\cO_{k',v})$-fixed vectors.

Let $\Pi$ be a cuspidal automorphic representation of $G_n(\bA_{k'})$ and let, for every place $v$ of $k$, $W_v$ be a Whittaker function in $\cW(\Pi_v,\Psi_{n,v})$ such that $W_v$ is $G_n(\cO_{k',v})$-invariant and satisfies $W_v(1)=1$ for almost all places $v$. Set
$$\displaystyle W=\prod_v W_v.$$
Then, $W$ is a well-defined function on $G_n(\bA_{k'})$ satisfying $W(ug)=\Psi_n(u)W(g)$ for every $(u,g)\in N_n(\bA_{k'})\times G_n(\bA_{k'})$. We define similarly $\widetilde{W}=\prod_v \widetilde{W}_v$. For every place $v$, choose a function $\phi_v\in \cS(k_v^n)$ such that $\phi_v=\mathbf{1}_{G_n(\cO_{k',v})}$ for almost all places $v$ and set $\phi=\prod_v \phi_v$ (a function on $\bA_{k'}^n$). We define similarly $\widehat{\phi}=\prod_v \widehat{\phi_v}$ where the local Fourier transforms are defined with respect to the local components $\Psi'_v$ of $\Psi'$. For $s\in \C$, we define, whenever convergent, the following global analog of the Zeta integrals discussed in the previous sections:
$$\displaystyle Z(s,W,\phi)=\int_{N_n(\bA)\backslash G_n(\bA)} W(h) \phi(e_n h) \lvert \det h\rvert^s dh.$$
Then, the following result is a consequence of \cite[Propositions 5 \& 6]{Kab}  (although, strictly speaking, the discussion of \textit{loc. cit.} assumes that $k'/k$ splits at all Archimedean places, this assumption is not used in the proof of the result below).

\begin{theo}\label{theo global functional eqn}
When $\Re(s)$ is sufficiently large, the integral defining $Z(s,W,\phi)$ is absolutely convergent and we have
$$\displaystyle Z(s,W,\phi)=\prod_v Z(s,W_v,\phi_v).$$
Moreover, the function $s\mapsto Z(s,W,\phi)$ has a meromorphic continuation to $\C$ that satisfies the functional equation
$$\displaystyle Z(s,W,\phi)=Z(1-s,\widetilde{W},\widehat{\phi}).$$
\end{theo}

\subsection{A globalization result}\label{Section globalization}

We continue with the setting fixed in Section \ref{Global functional equation}. Let $v_0,v_1$ be two distinct places of $k$ with $v_1$ non-Archimedean. Recall that in Section \ref{Section representations} we have equipped the set $\Temp(G_n(k'_{v_0}))$ of tempered representations of $G_n(k'_{v_0})$ with a topology. We have the following globalization result which is a simple consequence of the main result of \cite{FLM}. For convenience, we explain the deduction.

\begin{theo}\label{theo globalization}
Let $U$ be an open subset of $\Temp(G_n(k'_{v_0}))$. Then, there exists a cuspidal automorphic representation $\Pi$ of $G_n(\bA_{k'})$ such that $\Pi_{v_0}\in U$ and $\Pi_v$ is unramified for every non-Archimedean place $v\notin \{v_0,v_1 \}$.
\end{theo}

\noindent\ul{Proof}: Let $S_\infty$ be the set of all Archimedean places of $k$ and set $S=S_\infty\cup \{v_0 \}$, $k'_S=\prod_{v\mid S} k'_v$. Denote by $\Temp(G_n(k'_S))$ the set of all (isomorphism classes of) irreducible tempered representations of $G_n(k'_S)$ i.e. representations which can be written as (completed) tensor products of tempered irreducible representations of $G_n(k'_v)$ for every $v\in S$. Of course, we have an identification $\Temp(G_n(k'_S))=\prod_{v\in S} \Temp(G_n(k'_v))$ and we equip this set with the product topology. Clearly, it suffices to show that for every open $U$ of $\Temp(G_n(k'_S))$ there exists a cuspidal automorphic representation $\Pi$ of $G_n(\bA_{k'})$ such that $\Pi_S:=\widehat{\bigotimes_{v\in S}} \Pi_v$ belongs to $U$ and $\Pi_v$ is unramified for every $v\notin S\cup\{v_1 \}$. Let $G_n(\bA_{k'})^1$ and $G_n(k'_S)^1$ be the subgroups of elements $g$ in $G_n(\bA_{k'})$ and $G(k'_S)$ respectively satisfying $\lvert \det g\rvert_{\bA_{k'}}=1$ where $\lvert .\rvert_{\bA_{k'}}$ denotes the normalized absolute value on $\bA_{k'}$. Then, we have decompositions $G_n(\bA_{k'})=G_n(\bA_{k'})^1\times Z_\infty$ and $G_n(k'_S)=G_n(k'_S)^1\times Z_\infty$. Let $\Temp(G_n(k'_S)^1)$ be the set of isomorphism classes of tempered representations of $G_n(k'_S)^1$ i.e. representations obtained by restriction of representations in $\Temp(G_n(k'_S))$ that we equip with the quotient topology. Then, we have an isomorphism $\Temp(G_n(k'_S))\simeq \Temp(G_n(k'_S)^1)\times \widehat{Z_\infty}$. Let $U$ be an open subset of $\Temp(G_n(k'_S))$ that we may assume to be of the form $U_1\times V$ where $U_1\subset \Temp(G_n(k'_S)^1)$ and $V\subset \widehat{Z_\infty}$ are open. Let $\mathfrak{p}$ be the prime ideal in $\cO_{k'}$ corresponding to $v_1$ and for all $m\geqslant 1$ let $K^S(\mathfrak{p}^m)$ be the open subgroup $K_{v_1}(m)\times \prod_{v\notin S\cup\{v_1 \}} G_n(\cO_{k',v})$ of $\prod'_{v\notin S} G_n(k'_v)$ where $K_{v_1}(m):=\{k\in G_n(\cO_{k',v_1})\mid k\equiv I_n \mod{\mathfrak{p}^m} \}$. Let $\mu_{\Pl}$ be the Plancherel measure for $G_n(k'_S)^1$. Since $\mu_{\Pl}(U_1)$ is nonzero, by \cite[Theorem 2]{FLM} there exists an irreducible representation $\Pi'$ of $G_n(\bA_{k'})$ such that $\Pi'\mid_{G_n(\bA_{k'})^1}\hookrightarrow L^2(G_n(k')\backslash G_n(\bA_{k'})^1)$, $(\Pi'_S) \mid_{G_n(k'_S)^1}\in U_1$ and $(\Pi')^{K^S(\mathfrak{p}^m)}\neq 0$ for some $m\geqslant 1$. Let $\omega\in V$ and denote by $\Pi$ the unique irreducible representation of $G_n(\bA_{k'})$ whose restriction to $G_n(\bA_{k'})^1$ coincides with $\Pi'\mid_{G_n(\bA_{k'})^1}$ and whose central character equals $\omega$ on $Z_\infty$. Then $\Pi\hookrightarrow L^2(Z_\infty G_n(k')\backslash G_n(\bA_{k'}),\omega)$, $\Pi_S\in U$ and $\Pi^{K^S(\mathfrak{p}^m)}\neq 0$. Moreover, since $\Pi$ is tempered at every Archimedean place by \cite{Wall3} we have that $\Pi$ is actually a cuspidal (rather than just square-integrable) automorphic representation of $G_n(\bA_{k'})$. Obviously, $\Pi$ has all of the desired properties. $\blacksquare$

\begin{rem}
The proof of the above result uses the full strength of \cite{FLM} but in fact we can get a slightly weaker result which is however sufficient for the application we have in mind without appealing to the analysis of the continuous part of the spectral side of Arthur's weighted trace formula which is the hard part of \cite{FLM}. More precisely, the following ought to be provable (although the author hasn't checked all the details) using the limit property \cite[Proposition 3]{FLM} of the geometric side of Arthur's weighted trace formula together with Sauvageot's density principle \cite{Sau}:
\begin{center}
Let $G$ be any connected reductive group over $k$. Let $v_0$, $v_1$, $v_2$ be three distinct places of $k$ with $v_1$, $v_2$ non-Archimedean. Let $U$ be an open subset of $\Temp(G(k_{v_0}))$ and $\pi_1$ be a supercuspidal representation of $G(k_{v_1})$. Then, there exists a cuspidal automorphic representation $\Pi$ of $G(\bA)$ such that $\Pi_{v_0}\in U$, $\Pi_{v_1}$ is isomorphic to an unramified twist of $\pi_1$ and $\Pi_v$ is unramified for every non-Archimedean place $v\notin \{v_0,v_1 ,v_2\}$.
\end{center}
\end{rem}

\subsection{Proof of Theorem \ref{theo 1 inert} and Theorem \ref{theo 2 inert} in the nearly tempered case}\label{appendix proof of the theorem}

Let the setting be as in Section \ref{Section inert case}. We choose a quadratic extension $k'/k$ of number fields together with a place $v_0$ of $k$ such that:
\begin{itemize}
\item there exists an isomorphism $k'_{v_0}/k_{v_0}\simeq E/F$ that we fix from now on;
\item every Archimedean place of $k$ different from $v_0$ splits in $k'$.
\end{itemize}

Given a place $v$ of $k$, a nontrivial additive character $\psi'_v$ of $k_v$ and a generic irreducible representation $\pi_v$ of $G_n(k'_v)$ we let $\gamma^{Sh}(s,\pi_v,\As,\psi'_v)$ be the local Asai $\gamma$-factor defined by Shahidi (see \cite{Sha3}, \cite{Gold}). The goal of this section is to show the following:

\begin{theo}\label{theo comparison Shahidi}
For every $\pi\in \Irr_{\ntemp}(G_n(E))$, $W\in \cW(\pi,\psi_n)$, $\phi\in \cS(F^n)$ and $s\in \C$ with $\Re(s)=\frac{1}{2}$ we have
$$\displaystyle Z(1-s,\widetilde{W},\widehat{\phi})=\omega_\pi(\tau)^{n-1}\lvert \tau\rvert_E^{\frac{n(n-1)}{2}(s-1/2)}\lambda_{E/F}(\psi')^{-\frac{n(n-1)}{2}}\gamma^{Sh}(s,\pi,\As,\psi')Z(s,W,\phi).$$
\end{theo}

First we explain how to deduce from it Theorems \ref{theo 1 inert} and \ref{theo 2 inert} in the nearly tempered case. Let $\pi\in \Irr_{\ntemp}(G_n(E))$, $W\in \cW(\pi,\psi_n)$ and $\phi\in \cS(F^n)$. By the main results of \cite{Sha2} (in the Archimedean case) and \cite{Shank} (in the $p$-adic case) we have
\begin{align}\label{eq 1 proof tempered}
\displaystyle \gamma^{Sh}(s,\pi,\As,\psi')=\gamma(s,\pi,\As,\psi').
\end{align}
Therefore, as $s\mapsto \gamma(s,\pi,\As,\psi')$ is a meromorphic function on $\C$, the functional equation of Theorem \ref{theo comparison Shahidi} together with Lemma \ref{lem 2 conv Zeta integrals} implies that $s\mapsto Z(s,W,\phi)$ admits a meromorphic continuation to $\C$ and that the functional equation is satisfied for every $s\in \C$. This already shows Theorem \ref{theo 1 inert}. For Theorem \ref{theo 2 inert}, we rewrite the functional equation as
\begin{align}\label{eq 2 proof tempered}
\displaystyle \frac{Z(1-s,\widetilde{W},\widehat{\phi})}{L(1-s,\widetilde{\pi},\As)}=\omega_\pi(\tau)^{n-1}\lvert \tau\rvert_E^{\frac{n(n-1)}{2}(s-1/2)}\lambda_{E/F}(\psi')^{-\frac{n(n-1)}{2}}\epsilon(s,\pi,\As,\psi')\frac{Z(s,W,\phi)}{L(s,\pi,\As)}.
\end{align}
By Lemma \ref{lem 2 conv Zeta integrals} and since the functions $L(s,\pi,\As)$ and $L(1-s,\widetilde{\pi},\As)$ have no zeros, the functions $s\mapsto \frac{Z(s,W,\phi)}{L(s,\pi,\As)}$ and $s\mapsto \frac{Z(1-s,\widetilde{W},\widehat{\phi})}{L(1-s,\widetilde{\pi},\As)}$ are holomorphic on $\{\Re(s)\geqslant \frac{1}{2} \}$ and $\{\Re(s)\leqslant \frac{1}{2}\}$ respectively. Since the epsilon factor $\epsilon(s,\pi,\As,\psi')$ has no zeros either, this implies together with \ref{eq 2 proof tempered} that $s\mapsto \frac{Z(s,W,\phi)}{L(s,\pi,\As)}$ is holomorphic on $\{\Re(s)\geqslant \frac{1}{2} \}\cup \{\Re(s)\leqslant \frac{1}{2}\}=\C$. Moreover, as $\pi$ is nearly tempered the $L$-functions $L(s,\pi,\As)$ and $L(s,\widetilde{\pi},\As)$ are holomorphic in $\{\Re(s)\geqslant \frac{1}{2} \}$ (Lemma \ref{lem 2 As L functions}(ii)). Hence, by Lemma \ref{lem nonzero Zeta} for every $s_0\in \C$ with $\Re(s_0)\geqslant 1/2$ (resp. $\Re(s_0)\leqslant 1/2$) there exist $W\in \cW(\pi,\psi_n)$ and $\phi\in \cS(F^n)$ with $\frac{Z(s_0,W,\phi)}{L(s_0,\pi,\As)}\neq 0$ (resp. $\frac{Z(1-s_0,\widetilde{W},\widehat{\phi})}{L(1-s_0,\widetilde{\pi},\As)}\neq 0$). By \ref{eq 2 proof tempered} this implies that for every $s_0\in \C$ there exist $W\in \cW(\pi,\psi_n)$ and $\phi\in \cS(F^n)$ with $\frac{Z(s_0,W,\phi)}{L(s_0,\pi,\As)}\neq 0$. Finally, by Lemma \ref{lem 2 As L functions} (iii) and Lemma \ref{lem 2 conv Zeta integrals}, there exists $\epsilon>0$ such that both $s\mapsto \frac{Z(s,W,\phi)}{L(s,\pi,\As)}$ and $s\mapsto \frac{Z(s,\widetilde{W},\widehat{\phi})}{L(s,\widetilde{\pi},\As)}$ are of finite order in vertical strips of $\cH_{>\frac{1}{2}-\epsilon}$. Together with the functional equation \ref{eq 2 proof tempered}, this shows that $s\mapsto \frac{Z(s,W,\phi)}{L(s,\pi,\As)}$ is of finite order in vertical strips on all of $\C$. This ends the proof of Theorem \ref{theo 2 inert}.

We now proceed to the proof of Theorem \ref{theo comparison Shahidi}. Let $\Psi$ and $\Psi'$ be nontrivial characters of $k\backslash \bA$ and $k'\backslash \bA_{k'}$ with $\Psi$ trivial on $\bA$. By Lemma \ref{lem 1 inert} and the identity \ref{eq 1 proof tempered} we may assume that $\psi'=\Psi'_{v_0}$ and $\psi=\Psi_{v_0}$. We now show that it suffices to establish Theorem \ref{theo comparison Shahidi} for a dense subset $D\subseteq \Temp(G_n(E))$. Indeed, let $\pi\in \Irr_{\ntemp}(G_n(E))$ that we write as $\pi=i_{P(E)}^{G_n(E)}(\sigma_{\lambda_0})$ where $P=MU$ is a standard parabolic subgroup of $G_n$, $\sigma\in \Pi_2(M(E))$ and $\lambda_0=(\lambda_{0,1},\ldots,\lambda_{0,n})\in \cA_{M,\C}^*\subseteq \cA_\C^*=\C^n$ satisfies 
\begin{align}\label{eq 3 proof tempered}
\lvert \Re(\lambda_{0,i})\rvert<1/4,\;\;\; 1\leqslant i\leqslant n.
\end{align}
By Lemma \ref{lem unramified twists}, we may assume that $\lambda_0\in (\cA_{M,\C}^G)^*$. Set $\pi_\lambda=i_{P(E)}^{G_n(E)}(\sigma_\lambda)$ for every $\lambda\in \cA_{M,\C}^*$ (so that $\pi=\pi_{\lambda_0}$). Let $W\in \cW(\pi,\psi_n)$ and $\phi\in \cS(F^n)$. Let $\rho=(\frac{n-1}{2},\ldots,\frac{1-n}{2})\in (\cA^{G_n})^*$ be half the sum of the roots of $A_n$ in $B_n$. Then $\lvert \Re(\lambda_0)\rvert\prec \lvert \Re(\lambda_0)\rvert+\epsilon \rho$ for any $\epsilon>0$ and by \ref{eq 3 proof tempered} we can choose $\epsilon >0$ such that $-2\min(\lvert \Re(\lambda_0)\rvert+\epsilon \rho)<\frac{1}{2}$. Set $\mu=\lvert \Re(\lambda_0)\rvert+\epsilon \rho$ and $\cU[\prec \mu]=\{\lambda \in (\cA^{G_n}_{M,\C})^*\mid \lvert \Re(\lambda)\rvert\prec \mu\}$. By Corollary \ref{cor good sections Whittaker models}, there exists a map
$$\displaystyle \lambda\in (\cA_{M,\C}^G)^*\mapsto W_\lambda\in \cW(\pi_{\lambda},\psi_n)$$
such that $W_{\lambda_0}=W$, $W_\lambda\in \cC_\mu(N_n(E)\backslash G_n(E),\psi_n)$ for every $\lambda\in \cU[\prec \mu]$ and the induced map
$$\displaystyle \lambda\in \cU[\prec \mu]\mapsto W_\lambda\in \cC_\mu(N_n(E)\backslash G_n(E),\psi_n)$$
is analytic. Assume that Theorem \ref{theo comparison Shahidi} holds for a dense subset $D\subseteq \Temp(G_n(E))$. Then, since this theorem is insensitive to unramified twists (Lemma \ref{lem unramified twists}), there exists a dense subset $D_\sigma\subseteq i(\cA^{G_n}_M)^*$ such that Theorem \ref{theo comparison Shahidi} holds for $\pi_{\lambda}$ for every $\lambda\in D_\sigma$. In particular, we have
\begin{align}\label{eq 4 proof tempered}
\displaystyle Z(1-s,\widetilde{W}_\lambda,\widehat{\phi})=\omega_\pi(\tau)^{n-1}\lvert \tau\rvert_E^{\frac{n(n-1)}{2}(s-1/2)}\lambda_{E/F}(\psi')^{-\frac{n(n-1)}{2}}\gamma^{Sh}(s,\pi_\lambda,\As,\psi')Z(s,W_\lambda,\phi)
\end{align}
for every $\lambda\in D_\sigma$ and $s\in \C$ with $\Re(s)=\frac{1}{2}$. By the choice of $\mu$ and Lemma \ref{lem 1 conv Zeta integrals} for any $s\in \C$ with $\Re(s)=\frac{1}{2}$ the maps $\lambda\in \cU[\prec \mu]\mapsto Z(s,W_\lambda,\phi)$ and $\lambda\in \cU[\prec \mu]\mapsto Z(1-s,\widetilde{W}_\lambda,\widehat{\phi})$ are holomorphic. By \ref{eq 1 proof tempered} and Lemma \ref{lem 2 As L functions}(ii), the map $\lambda\in \cU[\prec \mu]\mapsto \gamma^{Sh}(s,\pi_\lambda,\As,\psi')$ is also holomorphic. Therefore, by density of $D_\sigma$ in $i(\cA_M^G)^*$, \ref{eq 4 proof tempered} is also satisfied for every $\lambda\in i(\cA_M^G)^*$ (by continuity) hence for every $\lambda\in \cU[\prec \mu]$ (by analyticity) and in particular for $\lambda=\lambda_0$. This shows that Theorem \ref{theo comparison Shahidi} holds for $\pi$ and ends the proof of the reduction.

Let $U\subset \Temp(G_n(E))$ be an open subset. We are now going to show that there exists $\pi\in U$ for which Theorem \ref{theo comparison Shahidi} holds. Let $v_1$ be a non-Archimedean place of $k$ which splits in $k'$. By Theorem \ref{theo globalization}, there exists a cuspidal automorphic representation $\Pi$ of $G_n(\bA_{k'})$ such that $\Pi_{v_0}\in U$ and $\Pi_v$ is unramified for every non-Archimedean place $v\notin \{v_0,v_1 \}$. Let $S=S_\infty\cup\{v_0,v_1 \}$ where $S_\infty$ stands for the set of all Archimedean places of $k$ and let $T$ be the (finite) set of places $v$ of $k$ outside $S$ which ramify in $k'$ or where $\Psi'_v$ or $\Psi_v$ is not unramified. Then, by \cite[Theorem 3.5(4)]{Sha3} the product
$$\displaystyle L^{S\cup T}(s,\Pi,\As)=\prod_{v\notin S\cup T} L(s,\Pi_v,\As),$$
which converges for $\Re(s)$ sufficiently large, admits a meromorphic continuation to $\C$ and we have
\begin{align}\label{eq 1 global identity}
\displaystyle L^{S\cup T}(s,\Pi,\As)=\prod_{v\in S\cup T} \gamma^{Sh}(s,\Pi_v,\As,\Psi'_v) L^{S\cup T}(1-s,\Pi^\vee,\As).
\end{align}
On the other hand, combining Theorem \ref{theo global functional eqn} with Lemma \ref{lem 1 unramified}, we find that
\begin{align}\label{eq 2 global identity}
\displaystyle \prod_{v\in S\cup T} Z(s,W_v,\phi_v) L^{S\cup T}(s,\Pi,\As)=\prod_{v\in S\cup T} Z(1-s,\widetilde{W}_v,\widehat{\phi}_v) L^{S\cup T}(1-s,\Pi^\vee,\As)
\end{align}
for every $(W_v)_{v}\in \prod_{v\in S\cup T}\cW(\Pi_v,\Psi_{n,v})$ and $(\phi_v)_v\in \prod_{v\in S\cup T}\cS(k_v^n)$. Notice that we haven't proven the meromorphic continuation of $Z(s,W_v,\phi_v)$ yet but, by Theorem \ref{theo global functional eqn} and the aforementioned result of Shahidi, we at least know that the products 
$$\prod_{v\in S\cup T} Z(s,W_v,\phi_v)\mbox{ and }\prod_{v\in S\cup T} Z(1-s,\widetilde{W}_v,\widehat{\phi}_v)$$
admit meromorphic continuations. For every place $v\in (S\cup T)\setminus \{ v_0\}$, either $v$ splits in $k'$ or $v$ is non-Archimedean and $\Pi_v$ is unramified. Therefore, by Theorem \ref{theo 1 split}, Corollary \ref{cor unramified} and Lemma \ref{lem nonzero Zeta} we can choose the local data so that the functions $s\mapsto Z(s,W_v,\phi_v)$ and $s\mapsto Z(1-s,\widetilde{W}_v,\widehat{\phi}_v)$ admit meromorphic continuation to $\C$ which are not identically zero and satisfying the functional equation
$$\displaystyle Z(1-s,\widetilde{W}_v,\widehat{\phi}_v)=\omega_{\Pi_v}(\tau)^{n-1}\lvert \tau\rvert_{k'_v}^{\frac{n(n-1)}{2}(s-1/2)}\lambda_{k'_v/k_v}(\Psi'_v)^{-\frac{n(n-1)}{2}}\gamma(s,\Pi_v,\As,\Psi'_v)Z(s,W_v,\phi_v)$$
for every $v\in (S\cup T)\setminus \{ v_0\}$. Here $\tau\in k'$ is the only trace zero element such that $\Psi(z)=\Psi'(\Tra_{\bA_{k'}/\bA}(\tau z))$ for all $z\in \bA_{k'}$. This already implies that the local Zeta integrals $Z(s,W,\phi)$ and $Z(1-s,\widetilde{W},\widehat{\phi})$ admit meromorphic continuation for every $W\in \cW(\Pi_{v_0},\psi_n)$ and $\phi\in \cS(F^n)$. Moreover, combining the above identities with the product formula
$$\displaystyle \prod_{v\in S\cup T} \omega_{\Pi_v}(\tau)^{n-1}\lvert \tau\rvert_{k'_v}^{\frac{n(n-1)}{2}(s-1/2)}\lambda_{k'_v/k_v}(\Psi'_v)^{-\frac{n(n-1)}{2}}=\prod_v \omega_{\Pi_v}(\tau)^{n-1}\lvert \tau\rvert_{k'_v}^{\frac{n(n-1)}{2}(s-1/2)}\lambda_{k'_v/k_v}(\Psi'_v)^{-\frac{n(n-1)}{2}}=1$$
(which is a consequence of the fact that $\Pi_v$, $k'_v/k_v$, $\Psi_v$ and $\Psi'_v$ are all unramified outside $S\cup T$) identity \ref{eq 2 global identity} becomes
\[\begin{aligned}
\displaystyle & \omega_{\Pi_{v_0}}(\tau)^{n-1}\lvert \tau\rvert_E^{\frac{n(n-1)}{2}(s-1/2)}\lambda_{E/F}(\psi')^{-\frac{n(n-1)}{2}}Z(s,W,\phi)L^{S\cup T}(s,\Pi,\As) \\
 & =Z(1-s,\widetilde{W},\widehat{\phi})\prod_{v\in S\cup T\setminus \{v_0 \}} \gamma(s,\Pi_v,\As,\Psi'_v) L^{S\cup T}(1-s,\Pi^\vee,\As)
\end{aligned}\]
for every $W\in \cW(\Pi_{v_0},\psi_n)$ and $\phi\in \cS(F^n)$. Taking the quotient with \ref{eq 1 global identity} and using \ref{eq 1 proof tempered}, we obtain
\[\begin{aligned}
\displaystyle & \omega_{\pi_{v_0}}(\tau)^{n-1}\lvert \tau\rvert_E^{\frac{n(n-1)}{2}(s-1/2)}\lambda_{E/F}(\psi')^{-\frac{n(n-1)}{2}}\gamma^{Sh}(s,\Pi_{v_0},\As,\psi')Z(s,W,\phi) \\
 & =Z(1-s,\widetilde{W},\widehat{\phi})\prod_{v\in S\cup T\setminus \{v_0 \}} \frac{\gamma(s,\Pi_v,\As,\Psi'_v)}{\gamma^{Sh}(s,\Pi_v,\As,\Psi'_v)}=Z(1-s,\widetilde{W},\widehat{\phi})
\end{aligned}\]
for every $W\in \cW(\Pi_{v_0},\psi_n)$ and $\phi\in \cS(F^n)$. This shows that Theorem \ref{theo comparison Shahidi} holds for $\Pi_{v_0}$ hence, as $U$ was arbitrary, for a dense subset of $\Temp(G_n(E))$. By the previous reduction, this proves Theorem \ref{theo comparison Shahidi} in all cases. $\blacksquare$

\subsection{Proof of Theorem \ref{theo 1 inert} and Theorem \ref{theo 2 inert} in the general case}\label{Section end of proof of main theorems}

Let $\pi$ be a generic irreducible representation of $G_n(E)$. There exists a standard parabolic subgroup $P=MU$ of $G_n$, $\sigma\in \Pi_2(M(E))$ and $\lambda_0\in \cA_{M,\C}^*$ such that $\pi=i_{P(E)}^{G_n(E)}(\sigma_{\lambda_0})$. By Lemma \ref{lem unramified twists}, we may assume that $\lambda_0\in (\cA^{G_n}_{M,\C})^*$. Set $\pi_\lambda=i_{P(E)}^{G_n(E)}(\sigma_{\lambda})$ for every $\lambda \in \cA_{M,\C}^*$. Let $W\in \cW(\pi,\psi_n)$ and $\phi\in \cS(F^n)$. Choose a section
$$\displaystyle \lambda\in (\cA_{M,\C}^G)^*\mapsto W_\lambda\in \cW(\pi_{\lambda},\psi_n)$$
as in Corollary \ref{cor good sections Whittaker models} with $W_{\lambda_0}=W$. For $\mu\in (\cA^{G_n})^*$ we define the open subset $\cU[\prec \mu]$ of $(\cA_{M,\C}^G)^*$ as in Corollary \ref{cor good sections Whittaker models}. In particular, for every $\mu\in (\cA^{G_n})^*$ the above section induces an analytic map
$$\displaystyle \lambda\in \cU[\prec \mu]\mapsto W_\lambda\in \cC_\mu(N_n(E)\backslash G_n(E),\psi_n).$$
Decomposing the $\epsilon$-factor $\epsilon(s,\pi_{\lambda},\As,\psi')$ according to Lemma \ref{lem 2 As L functions}(i) we see that there exists $C>0$, a linear form $L$ on $(\cA_{M,\C}^G)^*$ and an element $u\in \C$ of norm $1$ such that
$$\displaystyle \omega_\pi(\tau)^{n-1}\lvert \tau\rvert_E^{\frac{n(n-1)}{2}(s-1/2)}\lambda_{E/F}(\psi')^{-\frac{n(n-1)}{2}}\epsilon(s,\pi_\lambda,\As,\psi')=uC^{L(\lambda)+s-\frac{1}{2}}$$
for every $\lambda\in (\cA_{M,\C}^G)^*$ and $s\in \C$. Let $v$ be a square root of $u$ and set
$$\displaystyle \epsilon_{1/2}(s,\pi_{\lambda},\psi'):=v(\sqrt{C})^{L(\lambda)+s-\frac{1}{2}}, \;\;\; \lambda\in (\cA_{M,\C}^G)^*,s\in \C$$
so that by the above
\begin{align}\label{eq 1 proof general case}
\displaystyle \omega_\pi(\tau)^{n-1}\lvert \tau\rvert_E^{\frac{n(n-1)}{2}(s-1/2)}\lambda_{E/F}(\psi')^{-\frac{n(n-1)}{2}}\epsilon(s,\pi_\lambda,\As,\psi')=\epsilon_{1/2}(s,\pi_{\lambda},\psi')^2
\end{align}
for every $\lambda\in (\cA_{M,\C}^G)^*$ and $s\in \C$. We are now going to apply Proposition \ref{prop 1 holomorphic continuation} to the two following (for the moment partially defined) functions on $\C\times (\cA_{M,\C}^G)^*$:
$$\displaystyle Z_+(s,\lambda)=\epsilon_{1/2}(\frac{1}{2}+s,\pi_{\lambda},\psi')\frac{Z(\frac{1}{2}+s,W_\lambda,\phi)}{L(\frac{1}{2}+s,\pi_\lambda,\As)}$$
and
$$\displaystyle Z_-(s,\lambda)=\epsilon_{1/2}(\frac{1}{2}+s,\pi_{\lambda},\psi')^{-1}\frac{Z(\frac{1}{2}+s,\widetilde{W}_\lambda,\widehat{\phi})}{L(\frac{1}{2}+s,\widetilde{\pi_\lambda},\As)}.$$
Let $U\subseteq (\cA_{M,\C}^G)^*$ be the nonempty relatively compact connected open subset of vectors $\lambda=(\lambda_1,\ldots,\lambda_n)\in (\cA_{M,\C}^G)^*\subseteq \cA_\C^*$ satisfying $\lvert \lambda_i\rvert<\frac{1}{4}$ for every $1\leqslant i\leqslant n$. Then, for every $\lambda\in U$, the representation $\pi_\lambda$ is nearly tempered. Set as usual $\rho=(\frac{n-1}{2},\ldots,\frac{1-n}{2})\in (\cA^{G_n})^*$ for half the sum of the roots of $A_n$ in $B_n$. By Lemma \ref{lem 2 conv Zeta integrals} the functions $Z_+$, $Z_-$ are at least well-defined on $\{\Re(s)\geqslant \frac{1}{2} \}\times U$. Moreover, by the discussion of the previous section, $Z_+$ and $Z_-$ admit holomorphic extensions to $\C\times U$ which are of finite order in vertical strips in the first variable locally uniformly in the second variable and satisfying the functional equation
\begin{align}\label{eq 2 proof general case}
\displaystyle Z_+(s,\lambda)=Z_-(-s,\lambda),\;\; s\in \C, \lambda\in U.
\end{align}
Let $U'\subseteq (\cA_{M,\C}^G)^*$ be a relatively compact connected open subset containing $U$. Then, there exists $\mu\in (\cA^{G_n})^*$ such that $U'\subseteq \cU[\prec \mu]$. Hence, by Lemma \ref{lem 1 conv Zeta integrals} and Lemma \ref{lem 2 As L functions}(iii), there exists $C$ such that $Z_+$ and $Z_-$ admit holomorphic extensions to $\cH_{>C}\times U'$ which are of finite order in vertical strips in the first variable locally uniformly in the second variable. Therefore, $Z_+$ and $Z_-$ satisfy the hypothesis of Proposition \ref{prop 1 holomorphic continuation} and they extend to holomorphic functions on all of $\C\times (\cA_{M,\C}^G)^*$ of finite order in vertical strips in the first variable locally uniformly in the second variable satisfying \ref{eq 2 proof general case} for every $(s,\lambda)\in \C\times (\cA_{M,\C}^G)^*$. Specializing to $\lambda=\lambda_0$, by definition of $Z_+$, $Z_-$ and the relation \ref{eq 1 proof general case}, we deduce that $Z(s,W,\phi)$ and $Z(s,\widetilde{W},\widehat{\phi})$ admit meromorphic continuation to $\C$ satisfying the functional equation
$$\displaystyle  Z(1-s,\widetilde{W},\widehat{\phi})=\omega_\pi(\tau)^{n-1}\lvert \tau\rvert_E^{\frac{n(n-1)}{2}(s-1/2)}\lambda_{E/F}(\psi')^{-\frac{n(n-1)}{2}}\gamma(s,\pi,\As,\psi')Z(s,W,\phi).$$
This already proves Theorem \ref{theo 1 inert} for $\pi$. Moreover, we also obtain that $\displaystyle \frac{Z(s,W_\lambda,\phi)}{L(s,\pi_\lambda,\As)}$ is holomorphic of finite order in vertical strips thus implying Theorem \ref{theo 2 inert} for $\pi$.

\flushright Rapha\"el Beuzart-Plessis \\
Aix Marseille University \\
CNRS, Centrale Marseille, I2M\\
Marseille, France\\ 
email: rbeuzart@gmail.com

\end{document}